
\input amssym.def
\input amssym.tex


\def\item#1{\vskip1.3pt\hang\textindent {\rm #1}}

\def\itemitem#1{\vskip1.3pt\indent\hangindent2\parindent\textindent {\rm #1}}

\tolerance=300
\pretolerance=200
\hfuzz=1pt
\vfuzz=1pt


\hoffset=0.6in
\voffset=0.8in

\hsize=5.8 true in 
\vsize=8.5 true in
\parindent=25pt
\mathsurround=1pt
\parskip=1pt plus .25pt minus .25pt
\normallineskiplimit=.99pt

\countdef\revised=100
\mathchardef\emptyset="001F 
\chardef\ss="19
\def\3{\ss}
\def\anf{$\lower1.2ex\hbox{"}$}
\def\frac#1#2{{#1 \over #2}}
\def\>{>\!\!>}
\def\<{<\!\!<}

\def\into{\hookrightarrow}
\def\ssarr{\hbox to 30pt{\rightarrowfill}}
\def\sarr{\hbox to 40pt{\rightarrowfill}}
\def\arr{\hbox to 60pt{\rightarrowfill}}
\def\larr{\hbox to 60pt{\leftarrowfill}}
\def\Arr{\hbox to 80pt{\rightarrowfill}}

{}

\def\ad{\mathop{\rm ad}\nolimits}

\def\Ad{\mathop{\rm Ad}\nolimits}

\def\cone{\mathop{\rm cone}\nolimits}

\def\det{\mathop{\rm det}\nolimits}

\def\Exp{\mathop{\rm Exp}\nolimits}

\def\Hol{\mathop{\rm Hol}\nolimits}%
\def\id{\mathop{\rm id}\nolimits} 
\def\im{\mathop{\rm im}\nolimits}
\def\Im{\mathop{\rm Im}\nolimits}
\def\inf{\mathop{\rm inf}\nolimits}

\def\Int{\mathop{\rm int}\nolimits}

\def\Re{\mathop{\rm Re}\nolimits}

\def\Spec{\mathop{\rm Spec}\nolimits}

\def\sup{\mathop{\rm sup}\nolimits}

\def\tr{\mathop{\rm tr}\nolimits}

\def\0{{\bf 0}}
\def\1{{\bf 1}}

\def\a{{\frak a}}

\def\b{{\frak b}}
\def\cc{{\frak c}}

\def\e{{\frak e}}

\def\g{{\frak g}}

\def\h{{\frak h}}

\def\j{{\frak j}}
\def\k{{\frak k}}

\def\n{{\frak n}}

\def\p{{\frak p}}
\def\q{{\frak q}}

\def\sL{{\frak {sl}}}
\def\t{{\frak t}}

\def\z{{\frak z}}

\def\C{{\Bbb C}}

\def\N{{\Bbb N}}

\def\R{{\Bbb R}} 
 
\def\Z{{\Bbb Z}} 

\def\:{\colon}  
\def\.{{\cdot}}
\def\|{\Vert}
\def\bsk{\bigskip}

\def\giantskip{\vskip2\bigskipamount}
\def\gsk{\giantskip}
\def \la {\langle}
\def\msk{\medskip}
\def \ra {\rangle}
\def \res {\!\mid\!\!}

\def\ssk{\smallskip}

\def\bbr{\bigbreak}
\def\giantbreak{\par \ifdim\lastskip<2\bigskipamount \removelastskip
         \penalty-400 \giantskip\fi}

\def\nin{\noindent}
\def\cen{\centerline}
\def\pagebreak{\vskip 0pt plus 0.0001fil\break}
\def\linebreak{\break}

\def\hat{\widehat}

\def\epsilon{\varepsilon}
\def\eset{\emptyset}

\def\nin{\noindent}
\def\oline{\overline}

\def\pder#1,#2,#3 { {\partial #1 \over \partial #2}(#3)}
\def\pde#1,#2 { {\partial #1 \over \partial #2}}
\def\phi{\varphi}


\def\subeq{\subseteq}

\def\Rarrow{\Rightarrow}

\def\tilde{\widetilde}

\font\eightrm=cmr8


\font\smc=cmcsc10
\font\bfone=cmbx10 scaled\magstep1 
\font\bftwo=cmbx10 scaled\magstep2 

\def\qed{{\unskip\nobreak\hfil\penalty50\hskip .001pt \hbox{}\nobreak\hfil
          \vrule height 1.2ex width 1.1ex depth -.1ex
           \parfillskip=0pt\finalhyphendemerits=0\medbreak}\rm}

\def\qeddis{\eqno{\vrule height 1.2ex width 1.1ex depth -.1ex} $$
                   \medbreak\rm}

\def\Lemma #1. {\bigbreak\vskip-\parskip\noindent{\bf Lemma #1.}\quad\it}

\def\Sublemma #1. {\bigbreak\vskip-\parskip\noindent{\bf Sublemma #1.}\quad\it}

\def\Proposition #1. {\bigbreak\vskip-\parskip\noindent{\bf Proposition #1.}
\quad\it}

\def\Corollary #1. {\bigbreak\vskip-\parskip\nin{\bf Corollary #1.}
\quad\it}

\def\Theorem #1. {\bigbreak\vskip-\parskip\noindent{\bf Theorem #1.}
\quad\it}

\def\Definition #1. {\rm\bigbreak\vskip-\parskip\noindent{\bf Definition #1.}
\quad}

\def\Remark #1. {\rm\bigbreak\vskip-\parskip\noindent{\bf Remark #1.}\quad}

\def\Example #1. {\rm\bigbreak\vskip-\parskip\noindent{\bf Example #1.}\quad}

\def\Problems #1. {\bigbreak\vskip-\parskip\noindent{\bf Problems #1.}\quad}
\def\Problem #1. {\bigbreak\vskip-\parskip\noindent{\bf Problems #1.}\quad}

\def\Conjecture #1. {\bigbreak\vskip-\parskip\noindent{\bf Conjecture #1.}\quad}

\def\Proof#1.{\rm\par\ifdim\lastskip<\bigskipamount\removelastskip\fi\smallskip
            \noindent {\bf Proof.}\quad}

\def\Axiom #1. {\bigbreak\vskip-\parskip\noindent{\bf Axiom #1.}\quad\it}

\def\Satz #1. {\bigbreak\vskip-\parskip\noindent{\bf Satz #1.}\quad\it}

\def\Korollar #1. {\bbr\vskip-\parskip\nin{\bf Korollar #1.} \quad\it}

\def\Bemerkung #1. {\rm\bigbreak\vskip-\parskip\noindent{\bf Bemerkung #1.}
\quad}

\def\Beispiel #1. {\rm\bigbreak\vskip-\parskip\noindent{\bf Beispiel #1.}\quad}
\def\Aufgabe #1. {\rm\bigbreak\vskip-\parskip\noindent{\bf Aufgabe #1.}\quad}

\def\Beweis#1. {\rm\par\ifdim\lastskip<\bigskipamount\removelastskip\fi
           \smallskip\noindent {\bf Beweis.}\quad}

\nopagenumbers

\def\date{\ifcase\month\or January\or February \or March\or April\or May
\or June\or July\or August\or September\or October\or November
\or December\fi\space\number\day, \number\year}

\def\title{Title ??}
\def\author{Author ??}

\def\thanks#1{\footnote*{\eightrm#1}}

\def\rightheadline{\hfil{\eightrm\title}\hfil\tenbf\folio}
\def\leftheadline{\tenbf\folio\hfil{\eightrm\author}\hfil}
\headline={\vbox{\line{\ifodd\pageno\rightheadline\else\leftheadline\fi}}}

\def\firstheadline{}
\def\firstfootline{\cen{\rm\folio}}

\def\seite #1 {\pageno #1
               \headline={\ifnum\pageno=#1 \firstheadline
               \else\ifodd\pageno\rightheadline\else\leftheadline\fi\fi}
               \footline={\ifnum\pageno=#1 \firstfootline\else{}\fi}}

\newdimen\dimenone
 \def\checkleftspace#1#2#3#4{
 \dimenone=\pagetotal
 \advance\dimenone by -\pageshrink   
 \ifdim\dimenone>\pagegoal          
   \else\dimenone=\pagetotal
        \advance\dimenone by \pagestretch
        \ifdim\dimenone<\pagegoal
          \dimenone=\pagetotal
          \advance\dimenone by#1         
          \setbox0=\vbox{#2\parskip=0pt                
                     \hyphenpenalty=10000
                     \rightskip=0pt plus 5em
                     \noindent#3 \vskip#4}    
        \advance\dimenone by\ht0
        \advance\dimenone by 3\baselineskip   
        \ifdim\dimenone>\pagegoal\vfill\eject\fi
          \else\eject\fi\fi}


\def\subheadline #1{\nin\bigbreak\vskip-\lastskip
      \checkleftspace{0.7cm}{\bf}{#1}{\medskipamount}
          \indent\vskip0.7cm\centerline{\bf #1}\medskip}

\def\sectionheadline #1{\bigbreak\vskip-\lastskip
      \checkleftspace{1.1cm}{\bf}{#1}{\bigskipamount}
         \vbox{\vskip1.1cm}\cen{\bfone #1}\bsk}

\def\lsectionheadline #1 #2{\bigbreak\vskip-\lastskip
      \checkleftspace{1.1cm}{\bf}{#1}{\bigskipamount}
         \vbox{\vskip1.1cm}\cen{\bfone #1}\msk \cen{\bfone #2}\bsk}

\def\lchapterheadline #1 #2{\bigbreak\vskip-\lastskip\indent\vskip3cm
                       \cen{\bftwo #1} \msk \cen{\bftwo #2} \gsk}
\def\llsectionheadline #1 #2 #3{\bigbreak\vskip-\lastskip\indent\vskip1.8cm
\cen{\bfone #1} \msk \cen{\bfone #2} \msk \cen{\bfone #3} \nobreak\bsk\nobreak}


\newtoks\literat
\def\[#1 #2\par{\literat={#2\unskip.}%
\hbox{\vtop{\hsize=.15\hsize\nin [#1]\hfill}
\vtop{\hsize=.82\hsize\nin\the\literat}}\par
\vskip.3\baselineskip}

\mathchardef\emptyset="001F 
\def\address{Author: \tt$\backslash$def$\backslash$address$\{$??$\}$}

\def\firstpage{\nin
{\obeylines \parindent 0pt }
\vskip2cm
\centerline {\bfone \title}
\gsk
\centerline{\bf\author}

\vskip1.5cm \rm}

\def\addresstwo{}

\def\dlastpage{\par\vbox{\vskip1cm\nin
\line{
\vtop{\hsize=.5\hsize{\parindent=0pt\baselineskip=10pt\nin\address}}
\quad 
\vtop{\hsize=.42\hsize\nin{\parindent=0pt
\baselineskip=10pt\addresstwo}}
\hfill} }}


\def\title{Formal dimension for semisimple symmetric spaces}
\def\author{Bernhard Kr\"otz${}^*$}
\footnote{}{${}^*$Supported by the DFG-project HI 412/5-2}

\def\date{July 20, 1999}
\def\leftheadline{\tenbf\folio\hfil\eightrm\date}
\def\Box #1 { \msk\par\nin 
\centerline{
\vbox{\offinterlineskip
\hrule
\hbox{\vrule\strut\hskip1ex\hfil{\smc#1}\hfill\hskip1ex}
\hrule}\vrule}\msk }

\def\bs{\backslash} 

\def\address
{Bernhard Kr\"otz

Mathematisches Institut

Technische Universit\"at Clausthal

Erzstra\3e 1

D-38678 Clausthal-Zellerfeld

Germany

mabk@math.tu-clausthal.de}

\firstpage

\sectionheadline{Abstract}

If $G$ is a semisimple Lie group and $(\pi, {\cal H})$ an irreducible 
unitary representation of $G$ with square integrable matrix coefficients, 
then there exists a number $d(\pi)$ such that 
$$(\forall v, v',w, w'\in {\cal H})\qquad {1\over d(\pi)} \la v,v'\ra \la w', w\ra 
=\int_G \la \pi(g).v,w\ra \oline {\la\pi(g).v'.w'\ra } \ d\mu_G(g).$$ 
The constant $d(\pi)$ is called the {\it formal dimension} of $(\pi, {\cal H})$
and was computed by Harish-Chandra in [HC56, 66]. 
\par If now $H\bs G$ is a semisimple symmetric space and  $(\pi, {\cal H})$ an irreducible 
$H$-spherical unitary $(\pi, {\cal H})$ belonging to the holomorphic discrete series 
of $H\bs G$, then one can define a formal dimension $d(\pi)$ in an analogous manner. 
In this paper we compute $d(\pi)$ for these class of representations. 

\ssk\nin {\bf Keywords:} Holomorphic discrete series, highest weight representation, 
formal dimension, formal degree, spherical representation, $c$-functions. 
\ssk\nin {\bf AMS classification:} 22E46, 43A85.

\sectionheadline{Introduction}

Let $H\bs G$ be a semisimple irreducible simply connected non-compact symmetric
space admitting relative holomorphic discrete series, i.e., there
exists a unitary highest weight representation $(\pi_\lambda, {\cal
H}_\lambda)$ of $G$ and a non-zero $H$-invariant hyperfunction vector
$\nu\in {\cal H}_\lambda^{-\omega}$  such that 
$${1\over d(\lambda)}\:={1\over |\la \nu, v_\lambda\ra|^2 }\int_{HZ\bs
G} |\la \nu, \pi_\lambda(g).v_\lambda\ra|^2\ d\mu_{HZ\bs G} (HZg)$$
is finite. Here $v_\lambda$ denotes a highest weight vector, 
$Z$ the center of $G$ and $\mu_{HZ\bs G}$ a $G$-invariant measure 
on the homogeneous space $HZ\bs G$. Note that 
$v_\lambda$ and $\nu$ are unique up to scalar multiple as well as $\la \nu, v_\lambda\ra \neq 0$. 
Therefore the number $d(\lambda)$ is well defined and we call it the {\it
formal dimension} of the spherical highest weight 
representation $(\pi_\lambda, {\cal H}_\lambda)$.
We remark here that our definition of the formal dimension generalizes 
Harish-Chandra notion in the ``group case'', i.e., where $G=G_0\times G_0$ and 
$H=\Delta(G)=\{(g,g)\: g\in G_0\}$ for a simply connected hermitian 
Lie group $G_0$ (cf.\ [HC56] and Remark III.5 below). 
\par Note that 
the constants $d(\lambda)$ determine the part of the  Plancherel measure  
on $H\bs G$ which corresponds to the relative 
holomorphic discrete series. Thus the explicit knowledge of the 
formal dimensions gives us a better understanding of the Plancherel 
Theorem on $H\bs G$ which was recently obtained by van den Ban-Schlichtkrull
and Delorme (cf.\ [BS97,99], [De98]). 

\par Let $(\g,\tau)$ be the symmetric Lie algebra attached to $H\bs G$
and write $\g=\h\oplus\q$ for the $\tau$-eigenspace decomposition. 
If $\g=\k\oplus\p$ is a $\tau$-invariant Cartan decomposition 
of $\g$, then the algebraic characterization 
of $H\bs G$ admitting relative 
holomorphic discrete series is $\z(\k)\cap\q\neq 0$. 
Symmetric Lie algebras $(\g,\tau)$ having this property
are called {\it compactly causal} (cf.\
[Hi\'Ol96]). In the group case, i.e., $(\g,\tau)=(\g_0\oplus\g_0,\sigma)$
with $\sigma(X,Y)=(Y,X)$ the flip involution, this just means 
that $\g_0$ is hermitian. We remark here that the formal 
dimension in the group case was computed by Harish-Chandra
(cf.\ [HC56]).

\par In this paper we derive the formula for the formal dimension 
$d(\lambda)$ for compactly causal symmetric spaces. For the special class 
of Cayley type spaces this problem has been dealt with by Chadli
with Jordan algebra methods (cf.\ [Ch98]). The approach presented here is
general and purely representation theoretic. 

\par Our key result is the {\it Averaging Theorem} (cf.\ Theorem II.16)
which asserts that for large parameters $\lambda$ the $H$-integral 
over $v_\lambda$ converges.  More precisely, for large parameters $\lambda$
we prove that 
$$\int_H \pi_\lambda(h).v_\lambda\ d\mu_H(h)= {\la v_\lambda, v_\lambda\ra
\over \la\nu, v_\lambda\ra} c(\lambda+\rho) \nu,$$
where the left hand side has to be understood as a vector valued 
integral with values in the Fr\'echet space of 
hyperfunction vectors and 
$$c(\lambda)=\int_{\oline N\cap HAN} a_H(\oline n)^{-(\lambda+\rho)}\
d\mu_{\oline N}(\oline n)$$
denotes the $c$-function of the {\it non-compactly causal}
$c$-dual space $H^c\bs G^c$ (cf.\ [Hi\'Ol96]). 

\par To obtain the formula for the formal degree $d(\lambda)$, 
we plug in the relation for $\nu$ obtained from the Averaging Theorem in the
definition of $d(\lambda)$ and obtain for large parameters:

$$d(\lambda)=d(\lambda)^G c(\lambda+\rho),$$
where $d(\lambda)^G$ is the formal dimension of $(\pi_\lambda, {\cal
H}_\lambda)$ for the relative discrete series on $G$ (cf.\ Theorem
III.6).  Using some ideas  of \'Olafsson and \O rsted (cf.\
[\'O\O91]) we prove the analytic  continuation 
of our formula for 
$d(\lambda)$ (cf.\ Theorem IV.15).

\par The $c$-function $c(\lambda)$ can be written as a product 
$$c(\lambda)=c_0(\lambda) c_\Omega(\lambda),$$
where $c_0(\lambda)$ is the $c$-function of a certain Riemannian 
symmetric subspace of $H\bs G$ and $c_\Omega(\lambda)$  is the 
$c$-function of the real form $\Omega$ of the bounded 
symmetric domain ${\cal D}\cong G/K$. In particular we have 

$$d(\lambda)=d(\lambda)^G c_0(\lambda+\rho) c_\Omega(\lambda+\rho).$$
The ingredients in this formula of $d(\lambda)$ are known: Harish-Chandra computed $d(\lambda)^G$ in 
[HC56], Gindikin and Karpelevi\v c $c_0(\lambda)$ (cf.\ [GiKa62]) and finally \'Olafsson and the author 
computed $c_\Omega(\lambda)$ in [Kr\'Ol99] (see also [Fa95], [Gr97] for earlier results 
in important special cases). 

\par In the final section we give applications of our results 
to spherical holomorphic representation theory.  Recall that a unitary highest weight
representation $(\pi_\lambda, {\cal H}_\lambda)$ of $G$ extends naturally 
to a holomorphic representation of the maximal open complex 
Ol'shanski\u\i{} semigroup $S_{\rm max}^0=G\Exp(iW_{\rm max}^0)$ (cf.\
[Ne99b, Sect.\ XI.2]). If $(\pi_\lambda, {\cal H}_\lambda)$ is an
$H$-spherical unitary highest weight representation of $G$, then we 
define its {\it spherical character} by 
$$\Theta_\lambda\: S_{\rm max}^0\to\C, \ \ s\mapsto 
{\la v_\lambda, v_\lambda\ra\over |\la \nu, v_\lambda\ra|^2} 
\la \pi_\lambda(s).\nu, \nu\ra.$$ 
Note that $\Theta_\lambda$ is an $H$-biinvariant holomorphic 
function on $S_{\rm max}^0$. On the other hand on $S_{\rm max}^0\cap HAN$
one defines the {\it spherical function with parameter}
$\lambda\in\a_\C^*$
by 
$$\phi_\lambda\: S_{\rm max}^0\cap HAN\to \C,\ \  s\mapsto
\int_H a_H(sh)^{\lambda-\rho} \ d\mu_H(h),$$
whenever the right hand side makes
sense(cf.\ [FH\'O94] or [KN\'O98]). For large parameters $\lambda$ we prove the long searched relation 
of \'Olafsson (cf.\ [\'Ol98, Open Problem 7(1)])
$$(\forall s\in S_{\rm max}^0\cap HAN)\quad \Theta_\lambda(s)={1\over
c(\lambda+\rho)} \phi_{\lambda+\rho}(s)$$
(cf.\ Theorem V.4). 
Finally we want to point out that the results of this paper are a
major step towards a proof 
of the {\it Plancherel Theorem} of $G$-invariant Hilbert spaces of holomorphic 
functions on $G$-invariant subdomains of the Stein variety  $\Xi_{\rm
max}^0=G\times_H iW_{\rm max}^0$ (cf.\ [Ch98], 
[HiKr98, 99b], [H\'O\O 91], [Kr98, 99b], [KN\'O97], [Ne99a].)

\msk I am grateful to J. Faraut and J. Hilgert
who both made my stay in the stimulating atmosphere of Paris VI possible such that I could 
work on this problem.  Also I want to thank G. \'Olafsson 
for many exciting discussions on the subject and proofreading the manuscript.  
Finally I want to thank the referee for his very careful work.

\sectionheadline{I. Causal symmetric Lie algebras}

This subsection is a brief introduction to causal symmetric Lie
algebras. Purely algebraic definitions of ``causality'' are given and
the basic notation on the algebraic level is introduced.

\Definition I.1. Let $\g$ denote a finite dimensional Lie algebra 
over the real numbers. 
\par\nin (a) A {\it symmetric Lie algebra} is a pair $(\g,\tau)$,
where $\tau$ is an involutive automorphism of $\g$. We set 
$$\h\:=\{ X\in\g\: \tau(X)=X\}\quad\hbox{and}\quad \q\:=\{ X\in\g\:
\tau(X)=-X\}$$
and note that $\g=\h +\q$. We call $(\g,\tau)$ 
{\it irreducible}, if $\{0\}$ and $\g$ are the only $\tau$-invariant 
ideals of $\g$. 
\par\nin (b) We denote by $\g_\C$ the complexification of $\g$. If 
$\tau$ is a involution on  $\g$,  we also denote by $\tau$ the 
complex linear extension of $\tau$ to a endomorphism of $\g_\C$. 
\par\nin (c) The $c$-{\it dual} $\g^c$ of $(\g,\tau)$ is defined by 
$\g^c=\h +i\q$. 
\par\nin (d) If $\g$ is semisimple, then there exists a Cartan
involution $\theta$ of $\g$ which commutes with $\tau$ 
(cf.\ [Be57] or [KrNe96, Prop.\ I.5(iii)]). We write $\g=\k + \p$ 
for the corresponding Cartan decomposition. By subscripts we indicate 
intersections, for example $\h_\k=\h\cap\k$ etc. 
Since $\tau$ and $\theta $ commute,  we have  $\g=\h_\k + \h_\p +
\q_\k + \q_\p$. 
The prescription $\theta^c\:=\theta\tau\res_{\g^c}$ defines a Cartan
involution on $\g^c$  and we denote by $\g^c=\k^c+\p^c$ the
corresponding Cartan decomposition of $\g^c$. \qed

If $V\subeq \g$ is a subspace, then we set $\z(V)=\{X\in V\: (\forall
Y\in V)[X,Y]=0\}$.

\Definition I.2. Let $(\g,\tau)$ be an irreducible semisimple
symmetric Lie algebra and $\theta$ a Cartan involution of $\g$ commuting 
with $\tau$. Then we call $(\g,\tau)$

\item{(CC)} {\it compactly causal} if $\z(\q_\k)\neq\{0\}$.
\item{(NCC)} {\it non-compactly causal}, if $(\g^c,\tau\res_{\g^c})$ is (CC).
\item{(CT)} of {\it Cayley type}, if it is both (CC) and (NCC).\qed

\Lemma I.3. Let $(\g,\tau)$ be a symmetric Lie algebra. Then the
following assertions hold:

\item{(i)} The symmetric Lie algebra $(\g,\tau)$ is compactly
causal if and only if it  belongs to one of the following two
types:
\itemitem{(1)} The Lie algebra  $\g$ is simple hermitian and 
$\z(\k)\subeq \q$. 
\itemitem{(2)} The subalgebra $\h$ is simple hermitian and $(\g,\tau)\cong
(\h\oplus\h,\sigma)$, where $\sigma$ denotes the flip involution 
$\sigma(X,Y)=(Y,X)$. 

\item{(ii)} If $(\g,\tau)$ is compactly causal, then 
\itemitem{(a)} $\z(\k)\cap\q$ is one-dimensional,
\itemitem{(b)} every maximal abelian subspace $\b\subeq \q_\k$ is
maximal abelian in $\q$ and $\h_\p+\q_\k$.

\Proof. (i) This follows from [Hi\'Ol96, Lemma 1.3.5, Th.\ 1.3.8] or 
[KrNe96, Prop.\ V.6].
\par\nin (ii) This is a consequence of [Hi\'Ol96, Prop.\ 3.1.11].\qed

\Remark I.4. (a) From the view point of convex geometry and complex
analysis
the compactly causal symmetric spaces are the natural generalization of
hermitian groups in the symmetric space setting (cf.\ [Hi\'Ol96],
[KrNe96], [KN\'O97, 98] and [Ne99b]). The compactly and non-compactly
causal symmetric Lie algebras
have been classified; for a complete list see [Hi\'Ol96, Th.\ 3.2.8].
\par\nin (b) Suppose that $H\bs G$ is a simply connected 
symmetric space associated to 
an irreducible semisimple symmetric Lie algebra $(\g,\tau)$. 
If $(\g,\tau)$ is compactly causal, then Lemma I.3(ii)(b) implies
that the symmetric space $H\bs G$ admits relative holomorphic discrete
series (cf.\ [FJ80]). The converse is also true. This result seems
to us to be well known. But since we do not know a proof in the literature, 
we added a proof in Appendix B (cf.\ Lemma B.1).    
\qed

Let $(\g,\tau)$ be compactly causal. Recall
that this implies in particular that $\g$ is hermitian (cf.\ Lemma
I.3(i)). 

\par We choose
a maximal abelian subalgebra $i\a\subeq\q_\k$  and extend $i\a$ in $\k$ to
a compactly embedded Cartan subalgebra $\t$ of $\g$. Recall from Lemma
I.3(ii)(b) that $\a$ is maximal
abelian in $i\q$ and $\p^c$. Then 
$\t=\t_\h+ i\a$ and $\z(\k)\cap\q\subeq i\a$. By Lemma I.3 we know that 
$\z(\k)\cap\q=\R Z_0$ is one-dimensional and by [Hel78, Ch.\ VIII, \S 7]  
we can normalize $Z_0$ in such a way that $\Spec(\ad Z_0)=\{ -i, 0, i\}$ holds. 
We denote by $\hat\Delta=\hat\Delta(\g_\C, \t_\C)$ the root system of
$\g_\C$ with respect to $\t_\C$ and by $\Delta=\Delta(\g^c, \a)$ the 
restricted root system of $\g^c$ with respect to $\a$. Note that 
$\hat\Delta\res_\a\bs\{0\}=\Delta$. The corresponding root space 
decompositions are denoted by

$$\g_\C=\t_\C\oplus\bigoplus_{\hat\alpha\in \hat\Delta}
\g_\C^{\hat\alpha}\quad
\hbox{and} \quad \g^c=\a\oplus\z_\h(\a)\oplus \bigoplus_{\alpha\in \Delta}
(\g^c)^{\alpha}.$$
A root $\hat\alpha\in \hat\Delta$ is called {\it compact} if
$\hat\alpha(Z_0)=0$ and {\it non-compact} otherwise. Analogously one
defines compact and non-compact roots in $\Delta$. Write
$\hat\Delta_k$ and
$\hat\Delta_n$ for the set of all compact,
resp. non-compact, roots in $\hat \Delta$. Analogously one defines
$\Delta_k$ and $\Delta_n$.  

Once and for all we fix now a positive system $\hat\Delta^+$ of
$\hat\Delta$ such that
$$\hat\Delta_n^+\:=\hat\Delta^+\cap \hat\Delta_n=\{ \hat\alpha\in
\hat\Delta_n\: \hat\alpha(Z_0)=i\}$$
holds. A positive system $\Delta^+$ of $\Delta$ is defined by
$\Delta^+\:=\hat\Delta^+\res_\a\bs\{0\}$.

\sectionheadline{II. Spherical highest weight representations}

In this section we briefly recall the classification of analytic and 
hyperfunction vectors of a a unitary highest weight representation 
$(\pi_\lambda, {\cal H}_\lambda)$ of a simply connected compactly 
causal group $(G,\tau)$. Further we collect the basic facts of 
$H$-spherical highest weight representations. Then, after giving 
the definitions of the various $c$-functions associated to the 
non-compactly causal $c$-dual space $(G^c,\tau)$ of $(G,\tau)$, we
prove the key result of the whole paper: The Averaging Theorem, which
asserts that for large parameters $\lambda$ the $H$-integral over 
the highest weight vector converges in the Fr\'echet space of 
hyperfunction vectors. One obtains the up to 
scalar multiple uniquely determined $H$-spherical vector with a 
normalization constant which is given by a certain $c$-function.

\subheadline{Unitary highest weight representations}

Recall that if $G$ is a simply connected Lie
group associated to a symmetric Lie algebra $(\g,\tau)$, then $\tau$
integrates to an involution on $G$, also denoted by $\tau$, and
that the fixed point group $G^\tau$ is connected (cf.\ [Lo69, Th.\ 
3.4].) 

\par To a compactly causal symmetric Lie algebra 
$(\g,\tau)$ we associate the following analytic objects:

\msk
$$\matrix{ G & \hbox{simply connected Lie group with Lie
algebra $\g$},\hfill\cr 
G^c & \hbox{simply connected Lie group with Lie
algebra $\g^c$},\hfill\cr 
G_\C & \hbox{simply connected Lie group with Lie
algebra $\g_\C$},\hfill\cr
H & \hbox{$\tau$-fixed points in $G$},\hfill\cr 
H^c & \hbox{$\tau$-fixed points in $G^c$},\hfill\cr 
H_\C & \hbox{$\tau$-fixed points in $G_\C$},\hfill\cr 
K & \hbox{analytic subgroup in $G$ corresponding to $\k$}, \hfill\cr
K^c & \hbox{analytic subgroup in $G^c$ corresponding to $\k^c$}, \hfill\cr
K_\C & \hbox{analytic subgroup in $G_\C$ corresponding to $\k_\C$},\hfill\cr
H^0 & \hbox{centralizer of $\a$ in $H$},\hfill\cr
H^{c,0} & \hbox{centralizer of $\a$ in $H^c$},\hfill\cr
Z & \hbox{center of $G$}.\hfill\cr}$$
\msk 

\par\nin Note that even though both $H$ and $H^c$ are connected and
have the same Lie algebra, they are in general not isomorphic. Recall
that $Z\subeq K$.   
\par If $X$ is a topological space and $Y\subeq X$ is a subspace, then
we denote by $Y^0$ or $\Int Y$ the interior of $Y$ in $X$.  
\par For each $X\in \g_\C$ we denote by $\oline X$ the complex
conjugate of $X$ with respect to the real form $\g$.

\Definition II.1. (Complex Ol'shanski\u\i{} semigroups, cf.\ [Ne99, Ch.\ XI]) Let
$(\g,\tau)$ be a compactly causal symmetric Lie algebra and 
$\hat \Delta^+=\hat \Delta^+(\g_\C, \t_\C)$ be the positive 
system from Section I.
 
\par\nin (a) Associated to $\hat \Delta^+$ we define the {\it maximal cone} in $\t$ by 

$$\hat C_{\rm max}=\{ X\in \t\: (\forall \alpha\in \hat\Delta_n^+)\
i\alpha(X)\geq 0\}.$$
We set  $\hat W_{\rm max}\:=\oline{\Ad(G).\hat C_{\rm max}}$ and note that
$\hat W_{\rm max}$ is a closed convex $\Ad(G)$-invariant convex cone in
$\g$ admitting no affine lines and which is maximal with respect to 
these properties (cf.\ [Ne99b, Sect.\  VII.3]).  

\par\nin (b) Let $G_1\:=\la \exp_{G_\C}(\g)\ra$. By Lawson's Theorem 
$S_{\rm max,1}\:=G_1\exp_{G_\C}(i\hat W_{\rm max})$ is a closed subsemigroup of $G_\C$,
the {\it maximal complex Ol'shanski\u\i{} semigroup}, and the 
{\it polar map}
$$G_1\times \hat W_{\rm max}\to S_{\rm max,1}, \ \ (g, X)\mapsto g\exp(iX)$$
is a homeomorphism (cf.\ [La94, Th.\ 3.4]). 
\par Denote by $S_{\rm max}$ the universal covering
semigroup of $S_{\rm max,1}$ and write $\Exp\: \g+i\hat W_{\rm max}\to 
S_{\rm max}$ for the
lifting of $\exp_{G_\C}\res_{\g+i \hat W_{\rm max}}\: \g+i\hat W_{\rm max}\to
S_{\rm max,1}$. Then it is easy to see that $S_{\rm max}=G\Exp(i\hat W_{\rm max})$ and that
there is a polar map 
$$G\times \hat W_{\rm max}\to S, \ \ (g, X)\mapsto g\Exp(iX)$$
which is homeomorphism. We define the {\it interior} of
$S_{\rm max}$ by
$S_{\rm max}^0\:=G\Exp(i\hat W_{\rm max}^0)$. Note that $S_{\rm max}^0$ is an open semigroup
ideal in $S_{\rm max}$ which carries a natural complex structure for which
the semigroup multiplication is holomorphic. Further the prescription
$s=g\Exp(iX)\mapsto s^*=\Exp(iX) g^{-1}$ defines on $S_{\rm max}$ the structure 
of an {\it involutive semigroup}. Note that the involution is
antiholomorphic on $S_{\rm max}^0$.\qed

\Remark II.2.  Let $(\pi_\lambda, {\cal
H}_\lambda)$ be a unitary highest weight representation of $G$ with 
respect to $\hat\Delta^+$ and highest weight $\lambda\in
i\t^*$. Denote by $B({\cal H}_\lambda)$ the $C^*$-algebra of bounded operators
on ${\cal H}_\lambda$. Recall from [Ne99b, Th.\ XI.4.8] that $(\pi_\lambda,{\cal
H}_\lambda)$ has a  natural extension to a {\it holomorphic representation}
$\pi_\lambda\: S_{\rm max}\to B({\cal H}_\lambda)$ of $S_{\rm max}$,
i.e., $\pi_\lambda$ is strongly continuous, holomorphic when restricted
to $S_{\rm max}^0$ and satisfies $\pi_\lambda(s^*)=\pi_\lambda(s)^*$
for all $s\in S_{\rm max}$.
\par Note that for $X\in \hat W_{\rm max}$ one has 
$\pi_\lambda(\Exp(iX))=e^{id\pi_\lambda(X)}$.

\qed

\Definition II.3. Let $G$ be a Lie group and ${\cal H}$ a
Hilbert space. If $(\pi,{\cal H})$ is a unitary representation of $G$,
then we call  $v\in {\cal H}$ an {\it analytic vector} if the orbit
map $G\to {\cal H}, \ g\mapsto \pi(g).v$ is analytic. We denote by
${\cal H}^\omega$ the vector space of all analytic vectors of $(\pi,
{\cal H})$. There is a natural locally convex topology
on ${\cal H}^\omega$ for which the representation $(\pi^{\omega},
{\cal H}^\omega)$ of $G$ on ${\cal H}^\omega$ is continuous (cf.\ [KN\'O97, Appendix]). The
strong antidual of ${\cal H}^\omega$ is denoted by 
${\cal H}^{-\omega}$ and the elements of ${\cal H}^{-\omega}$ are called {\it hyperfunction
vectors}. Note that there is a 
natural chain of continuous inclusions
$${\cal H}^\omega\into {\cal H}\into {\cal H}^{-\omega}.$$
The natural extension of $(\pi,{\cal H})$ to a representation on the space 
of hyperfunction vectors is denoted by $(\pi^{-\omega}, {\cal
H}^{-\omega})$ and given explicitly by 
$$ \la \pi^{-\omega}(g).\nu, v \ra := \la \nu, \pi^\omega(g^{-1}).v
\ra. \qeddis

\Proposition II.4. Let $(\pi_\lambda, {\cal H}_\lambda)$ be a unitary 
highest weight representation of $G$ with respect to $\hat\Delta^+$
and highest weight $\lambda$. Let $X\in \Int \hat W_{\rm max}$ be an
arbitrary element. Then the analytic vectors of $(\pi_\lambda, {\cal
H}_\lambda)$ are given by 

$${\cal H}_\lambda^\omega=\bigcup_{t>0} \pi_\lambda(\Exp(tiX)).{\cal
H}_\lambda$$
and the topology on ${\cal H}_\lambda^\omega$ is the finest locally
convex topology making for all $t>0$ the maps
${\cal H}_\lambda\to {\cal H}_\lambda^\omega, \ v\mapsto \pi_\lambda(\Exp(tiX)).v$
continuous.

\Proof. [KN\'O98, Prop.\ A.5].\qed

If $\lambda\in i\t^*$ is  dominant integral for $\hat \Delta_k^+$, we 
denote by $(\pi_\lambda^K, F(\lambda))$ the irreducible highest weight
representation of $K$ with highest weight $\lambda$.   Note that
$(\pi_\lambda^K, F(\lambda))$ extends naturally to a holomorphic
representation of the universal covering group  $\tilde{K_\C}$ of
$K_\C$ and which we denote by the same symbol.

\Remark II.5.  Recall that $(\pi_\lambda, {\cal H}_\lambda)$ can be  
realized in the Fr\'echet space $\Hol({\cal D}, F(\lambda))$ of 
$F(\lambda)$-valued holomorphic 
functions on the Harish-Chandra realization ${\cal D}$ of the
hermitian symmetric space $G/K$. So let us assume ${\cal H}_\lambda\subeq 
\Hol({\cal D}, F(\lambda))$. Then for all $z\in {\cal D}$ and
$v\in F(\lambda)$ the point evaluation 

$${\cal H}_\lambda\to \C, \ \ f\mapsto \la f(z),
v\ra $$
is continuous, hence can be written as $\la f(z),v\ra =\la f,
K_{z,v}^\lambda\ra $ for some $K_{z,v}^\lambda\in {\cal H}_\lambda$. One
can show that all vectors $K_{z,v}^\lambda$ are analytic. 
Then, if $\oline{{\cal H}_\lambda}$ denotes the closure of ${\cal H}_\lambda$ in the
nuclear Fr\'echet space $\Hol({\cal D}, F(\lambda))$, then the mapping

$$r\: {\cal H}_\lambda^{-\omega}\to \Hol({\cal D},
F(\lambda)), \ \ \nu\mapsto r(\nu);\ \la r(\nu)(z), v\ra=\nu(K_{z,v}^\lambda)$$
is a $G$-equivariant topological isomorphism onto its closed image
$\im r= \oline{{\cal H}_\lambda}$. In particular ${\cal
H}_\lambda^{-\omega}$ is a nuclear Fr\'echet space (cf.\ [Kr99a, 
Sect.\ III] for all that).\qed

\subheadline{Spherical representations}

\Definition II.6. Let $G$ be a Lie group, $H\subeq G$ a closed subgroup
and $(\pi,{\cal H})$ a unitary representation of $G$. Then we 
write $({\cal H}^{-\omega})^H$ for the set of 
all those elements $\nu\in {\cal H}^{-\omega}$ satisfying $\pi^{-\omega}(h).\nu=\nu$
for all $h\in H$. The unitary representation $(\pi,{\cal H})$ is called 
$H$-{\it spherical} if there exists a cyclic vector 
$\nu\in ({\cal H}^{-\omega})^H$ for $(\pi^{-\omega}, {\cal H}^{-\omega})$.\qed

For $\lambda\in i\t^*$  dominant integral with respect to $\hat \Delta_k^+$
recall the definition of the {\it generalized Verma module}

$$N(\lambda)\:={\cal U}(\g_\C)\otimes_{{\cal U}(\k_\C\ltimes \p^+)} F(\lambda)$$ 
which is a highest weight module of $\g$ with respect to
$\hat\Delta^+$ and highest weight $\lambda$
(cf.\ [EHW83]). We denote by $L(\lambda)$ the unique
irreducible quotient of $N(\lambda)$.

\Proposition II.7. Let $(\pi_\lambda, {\cal H}_\lambda)$ be a unitary
highest weight representation of $G$ with respect to  $\hat\Delta^+$. 

\item{(i)} If $(\pi_\lambda, {\cal H}_\lambda)$ is $H$-spherical, then
$(\pi_\lambda^K, F(\lambda))$ is $H\cap K$-spherical. In particular
$\lambda\in \a^*$ and the highest weight vector $v_\lambda\in {\cal
H}_\lambda$ is fixed by $H^0$. 
\item{(ii)} The restriction mapping 
$${\rm Res}\: ({\cal H}^{-\omega})^H\to F(\lambda)^{H\cap K}, \ \ \nu\mapsto\nu\res_{F(\lambda)} $$
is injective. In particular, $\dim({\cal H}^{-\omega})^H\leq 1$ and $\la \nu, v_\lambda\ra \neq 0$
for $\nu\neq 0$. 
Moreover if $L(\lambda)=N(\lambda)$, then ${\rm Res}$ is a bijection, i.e.,
$(\pi_\lambda, {\cal H}_\lambda)$ is $H$-spherical if and only if 
$(\pi_\lambda^K, F(\lambda))$ is $H\cap K$-spherical.

\Proof. (i) is a special case of [KN\'O97, Prop.\ VI.5] and
(ii) a special case of [Kr99a, Th.\ III.14].\qed

\Remark II.8. In general it is not true that $(\pi_\lambda, {\cal
H}_\lambda)$ is $H$-spherical if the minimal $K$-type $(\pi_\lambda^K,
F(\lambda))$ is $H\cap K$-spherical. For a counter example see [Kr99a,
Ex.\ III.16].\qed

\subheadline{ The $c$-functions on the $c$-dual space $H^c\bs G^c$.}

To the positve system $\Delta^+=\Delta^+(\g^c, \a)$ we associate several subalgebras of $\g^c$
$$\matrix{\n=\bigoplus_{\alpha\in \Delta^+} (\g^c)^\alpha,&\qquad  
\oline \n=\bigoplus_{\alpha\in \Delta^-} (\g^c)^\alpha,\cr
& \cr
\n_n^{\pm}=\bigoplus_{\alpha\in \Delta_n^{\pm}} (\g^c)^\alpha,
&\qquad \n_k^\pm=\bigoplus_{\alpha\in \Delta_k^\pm} (\g^c)^\alpha.\cr}$$  
Further we set 
$$\p^\pm\:=\bigoplus_{\hat\alpha\in \hat\Delta_n^+}\g_\C^{\hat\alpha}
\quad \hbox{and}\quad \g(0)\:=\h_\k + i\q_\k\subeq \g^c.$$

\Remark II.9. (a) The subalgebras $\p^\pm$ and $\k_\C$ of
$\g_\C$ are invariant under complex conjugation with respect to $\g^c$
and we have $\p^\pm\cap\g^c=\n_n^\pm$ as well as
$\k_\C\cap\g^c=\g(0)$. Thus the decomposition
$\g_\C=\p^+\oplus\k_\C\oplus\p^-$
induces a splitting in subalgebras of $\g^c$ 
$$\g^c=\n_n^+\oplus \g(0)\oplus \n_n^-.$$
\par\nin (b) Recall that $\g^c=\h\oplus\a\oplus\n$. The $\h\oplus\a\oplus\n$-
decomposition restricted to $\g(0)$ coincides with an Iwasawa
decomposition of $\g(0)$ given by $\g(0)=\k(0)\oplus \a \oplus
\n_k^+$, where $\k(0)\:=\h\cap\g(0)=\k^c\cap\g(0)$.  \qed

We let $H_\C\cap K_\C$ act on $H_\C\times K_\C$ from the left 
by $x.(h, k)\:=(hx^{-1}, xk)$ and denote by 
$$M\:=H_\C\times_{H_\C\cap K_\C} K_\C$$
the corresponding quotient space. The $H_\C\cap K_\C$-coset of an element
$(h,k)\in K_\C\times H_\C$ is denoted by $[h,k]$. 
If $\tilde{H_\C}$ and $\tilde{K_\C}$ denote the universal coverings of
$H_\C$ and $K_\C$, respectively, then we realize the universal cover
$\tilde M$ of $M$ by 
$$\tilde M =\tilde{H_\C}\times_{(\tilde{H_\C}\cap \tilde{K_\C})_0}
\tilde {K_\C}.$$
Further let $P^\pm\:=\exp_{G_\C}(\p^\pm)$. Recall that $\p^\pm$ are abelian and
that the exponential mapping $\exp_{G_\C}\res_{\p^\pm}\: \p^\pm\to P^\pm$ is  an
isomorphism. In particular $P^\pm$ is simply connected.

\Proposition II.10. {\rm (The $H_\C K_\C P^+$-decomposition)} The
following assertions hold:

\item{(i)} The multiplication mapping 
$$ M\times P^+\to G_\C,\ \ ([h, k], p_+)\mapsto hkp_+$$
is a biholomorphic map onto its open image $H_\C K_\C
P^+$. Furthermore:
\itemitem{(a)} The open submanifold $H_\C K_\C P^+$ is dense in $G_\C$
with complement of Haar measure zero. 
\itemitem{(b)} We have $S_{\rm max,1}\subeq H_\C
K_\C P^+$.
\item{(ii)} If $j\: S_{\rm max,1}\to M\times P^+$ denotes the injection
obtained from the isomorphism in {\rm (i)}, then $j$ lifts to an
inclusion mapping $\tilde j\: S_{\rm max}\to \tilde M\times P^+$.

\Proof. (i) [KN\'O97, Prop.\ II.6, Lemma III.7]. 
\par\nin (ii) Since $\pi_1(S_{\rm max,1})=\pi_1(G_1)\subeq Z(G)\subeq Z(K)$, it suffices
to show that $\tilde j\res_K$ is injective. We may assume that  $K\subeq
\tilde{K_\C}$, since both $K$ and $\tilde {K_\C}$ are simply connected
and $\k$ is a maximal compact subalgebra of $\k_\C$. Further $K$
normalizes $P^+$, and so  establishing the injectivity of $\tilde
j\res_K$ boils down to proving injectivity of  $K\to \tilde M, \
k\mapsto [\1, k]$, which is obvious. \qed

\par We denote by $G(0)$, $A$, $N$, $\oline N$, $N_k^\pm$ and
$N_n^\pm$ the analytic subgroups of $G^c$ corresponding to 
$\g(0)$, $\a$, $\n$, $\oline \n$, $\n_k^\pm$ and
$\n_n^\pm$.

\Remark II.11. (a) In view of the Bruhat 
decomposition of $\tilde {K_\C}$, we may identify $AN_k^+$ as a
subgroup of $\tilde {K_\C}$. Note that $N=N_k^+\rtimes N_n^+$ and so
every $n\in N$ can be written uniquely as $n=n_k n_n$ with $n_k\in N_k^+$ and
$n_n\in N_n^+$.
Thus we conclude from Proposition II.10(ii) that the map
$$H\times A\times N\to \tilde M\times P^+, \ \ (h,a, n_kn_n)\mapsto 
([h,an_k], n_n)$$
is an analytic diffeormphism onto its image which we denote by
$HAN$. Accordingly
every element $s\in HAN$ can be written uniquely as $s=h_H(s) a_H(s)
n_H(s)$ with $h_H(s)\in H$, $a_H(s)\in A$ and $n_H(s)\in N$ all
depending analytically on $s\in HAN$. 
\par\nin (b) If ${\cal D}\subeq \p^+$ denotes the Harish-Chandra 
realization of the hermitian symmetric space $G/K$  and $\oline {\cal D}$ its
conjugate in $\p^-$, then we set $\Omega\:=\oline{\cal D}\cap \n_n^-$.
In the sequel we realize $\Omega$ as a subset of $N_n^-$ via the 
exponential mapping. Recall from [KN\'O98, Lemma I.18] that 
$$H^cAN=\Omega G(0) N_n^+\quad\hbox{and}\quad \oline N\cap
H^cAN=\Omega\rtimes N_k^-.$$
On the other hand $\Omega$ can also naturally be realized in 
$\tilde M\times P^+$. In particular we 
obtain that the submanifold $\Omega\rtimes N_k^-$ of $\oline N$ is
naturally included in  $\tilde M\times P^+$. Denote this realization by
$\oline N\cap HAN$. Further the $HAN$-decomposition and the 
$H^c AN$-decomposition (cf.\ [KN\'O97, Prop.\ II.4]) coincide on 
$\oline N\cap HAN$. In the sequel we will use this fact frequently 
without mentioning it. 
\par\nin (c) Let $p\: X\to H^cAN$ the universal covering of
$H^cAN$. Since $X$ is simply connected, there exists a natural regular
map $\pi\: X\to \tilde M\times P^+$ with $\pi(X)=HAN$. In particular, 
the prescription 
$$K^c\cap HAN\:=\pi(p^{-1}(K^c\cap H^cAN))$$
defines an open submanifold of $HAN$. \qed

Note that the exponential mapping
$\exp_{\tilde {K_\C}}\res_\a\: \a\to A$ is an isomorphism, hence has an
inverse which we denote by $\log\: A\to \a$. For each $\lambda\in
\a_\C^*$ and $a\in A$ we set $a^\lambda=e^{\lambda(\log a)}$.

\Definition II.12. (The $c$-functions) We write $\rho$, $\rho_k$ and $\rho_n$ 
for the elements of $\a^*$ given  by 
${1\over 2} \tr \ad_\n$, ${1\over 2} \tr \ad_{\n_k^+}$ and 
${1\over 2} \tr \ad_{\n_n^+}$, respectively.
To $\lambda\in \a_\C^*$ we associate the following $c$-{\it functions}: 

$$\eqalign{c(\lambda)&\:=\int_{\oline N\cap (HAN)} a_H(\oline n)^{-(\lambda+\rho)}
\ d\mu_{\oline N} (\oline n),\cr 
c_{\Omega}(\lambda)&\:=\int_{\Omega} a_H(\oline n)^{-(\lambda+\rho)}
\ d\mu_{N^-_n} (\oline n), \cr}$$
and
$$c_0(\lambda)\:=\int_{N_k^-} a_H(\oline n)^{-(\lambda+\rho_k)}
\ d\mu_{N_k^-} (\oline n)$$
provided the integrals converge absolutely (cf.\ [FH\'O94] and  [KN\'O98]). 
We write ${\cal E}$ for the set of all $\lambda\in\a_\C^*$ for which the defining integral for $c$ 
converges absolutely.  Accordingly we define ${\cal E}_\Omega$  
and ${\cal E}_0$. Note that $c_0$ is the $c$-function of the
non-compact Riemannian symmetric space $K(0)\bs G(0)$, where 
$K(0)\:=G(0)^\tau$. \qed

For each $\alpha\in \Delta^+$  let $\check\alpha\in \a$ be the 
corresponding {\it coroot}, i.e., $\check\alpha\in [(\g^c)^\alpha, 
(\g^c)^{\tau\alpha}]$ such that $\alpha(\check\alpha)=2$. 
Associated to $\Delta^+$ we define two minimal cones in $\a$ by 
$$C_{\rm min}\:=\cone(\{ \check\alpha\: \alpha\in \Delta_n^+\})\quad 
\hbox{and}\quad \check C_k\:=\cone(\{\check\alpha\: \alpha\in \Delta_k^+\}).$$

\Definition II.13. Let $V$ be a finite dimensional vector space and
$V^*$ its dual.

\par\nin (a) If $C\subeq V$ is a convex set, then its {\it limit cone}
is defined by $\lim C=\{ x\in V\: x+C\subeq C\}$. Note that $\lim C$
is a convex cone and that $\lim C$ is closed if $C$ is open or
closed. 
\par\nin (b) If $E\subeq V$ is a subset, then its {\it dual cone}
is defined by $E^\star\:=\{ \alpha\in V^*\: \alpha\res_V\geq
0\}$. Note that $E^\star$ is a closed convex cone in $V^*$. \qed

\Theorem II.14. The various $c$-functions are related by 
$$c(\lambda)=c_0(\lambda) c_\Omega(\lambda)$$
and ${\cal E}={\cal E}_\Omega\cap {\cal E}_0$. Further:
\item{(i)} The domain of convergence ${\cal E}_\Omega$ of $c_\Omega$
is a tube domain ${\cal E}_\Omega=i\a^*+ {\cal E}_{\Omega, \R}$
with 
$${\cal E}_{\Omega,\R}=\{ \lambda\in \a^*\: (\forall \alpha\in 
\Delta_n^+)\  \lambda(\check\alpha)<2-m_\alpha\}, $$
where $m_\alpha\:=\dim (\g^c)^\alpha$. Further for all 
$\lambda\in {\cal E}_\Omega$ we have 
$$c_\Omega(\lambda)=\prod_{\alpha\in \Delta_n^+} B\Big(-{\lambda(\check\alpha)\over 2} -
{m_\alpha\over 2} +1, {m_\alpha\over 2}\Big), $$
where $B$ denotes the Euler Beta-function. In particular: 
\itemitem{(a)} $-\rho-C_{\rm min}^\star\subeq {\cal E}_{\Omega, \R}$
and $\lim {\cal E}_{\Omega, \R}=-C_{\rm min}^\star$.
\itemitem{(b)} The function $c_\Omega$ is holomorphic on ${\cal
E}_\Omega$ and $c_\Omega\res_{{\cal E}_\Omega +\mu}$ is bounded for all
$\mu\in -\rho-C_{\rm min}^\star$. 
\item{(ii)} The domain of convergence of $c_0$ is given by 
$${\cal E}_0=i\a^* +\Int \check C_k^\star,$$
$c_0$ is holomorphic on ${\cal E}_0$ and $c_0\res_{{\cal
E}_0+\mu}$ is bounded for all $\mu\in \rho_k+\check C_k^\star$.

\Proof. The  product formula $c(\lambda)=c_0(\lambda)
c_\Omega(\lambda)$ and the relation ${\cal E}={\cal E}_\Omega\cap {\cal E}_0$ 
are a special case of [KN\'O98, Lemma IV.5]. 
\par\nin (i) [Kr\'Ol99, Th.\ III.5]. 
\par\nin (ii) All this follows from the Gindikin-Karpelevic product 
formula for $c_0$  (cf.\  [Hel84, \ Ch.\ IV, Th.\ 6.13]).\qed

\subheadline{The Averaging Theorem}

\Lemma II.15. The group $H^0$ is compact and up to normalization of Haar measures for all  
$f\in L^1(H/H^0)$ the following integration
formulas hold:
\item{(i)}$$\int_H f(hH^0)\ d\mu_H(h)=
\int_{\oline N\cap (HAN)} f(h_H(\oline n)H^0) a_H(\oline n)^{-2\rho}\ d\mu_{\oline N} (\oline n).$$
\item{(ii)} $$\int_H f(hH^0)\ d\mu_H(h)=\int_{K^c\cap (HAN)} 
f(h_H(k)H^0) a_H(k)^{-2\rho}\ d\mu_{K^c}(k).$$

\Proof. In [KN\'O98, Lemma III.15(i)] it is proved that $H^{c,0}$
is compact and exactly the same argument also yields that
$H^0$ is compact. In view of this fact and our identifications of the
various domains in the big complex manifold $\tilde M\times P^+$
(cf.\ Remark II.11),
(i) follows from [KN\'O98, Prop.\ 1.19] and (ii) from [\'Ol87, Lemma 1.3].\qed

\Theorem II.16. {\rm (The Averaging Theorem)} 
Let $(\pi_\lambda, {\cal H}_\lambda)$ be a unitary
highest weight representation of $G$ for which $(\pi_\lambda^K,
F(\lambda))$ is $H\cap K$-spherical. If $v_\lambda$ is  a highest
weight vector, then the  vector valued integral 
$$\int_H \pi_\lambda(h).v_\lambda\ d\mu_H(h)$$
with values in the Fr\'echet space ${\cal H}_\lambda^{-\omega}$ (cf.\
{\rm Remark II.5}) converges and defines a non-zero $H$-fixed
hyperfunction vector if and only if $\lambda+\rho \in {\cal
E}_\Omega$. If this condition is satisfied and  $0\neq \nu \in({\cal
H}_\lambda^{-\omega})^H$, then 

$$\int_H \pi_\lambda(h).v_\lambda\ d\mu_H(h)= {\la v_\lambda, v_\lambda\ra
\over \la\nu, v_\lambda\ra} c(\lambda+\rho) \nu.$$

\Proof. Step 1: The analytic function $S_{\rm max}^0\cap HAN\to {\cal
H}_\lambda, \ s\mapsto \pi_\lambda(s).v_\lambda$ extends to an
analytic function $F\: HAN\to {\cal H}_\lambda$ and is given
ecplicitly  by 
$F(s)=a_H(s)^\lambda \pi_\lambda(h_H(s)).v_\lambda$.
\par In fact since $d\pi_\lambda(X).v_\lambda =0$ for all $X\in\n$,
the standard argument of differentiating yields
$$\pi_\lambda(s).v_\lambda=\pi_\lambda(h_H(s)a_H(s)n_H(s)).v_\lambda
=\pi_\lambda(h_H(s)a_H(s)).v_\lambda=a_H(s)^\lambda\pi_\lambda(h_H(s)).v_\lambda,$$
establishing Step 1.

\par\nin Step 2: The integral exists if and only if $\lambda +\rho\in
{\cal E}_\Omega$.
\par   Let $X\in \Int \hat W_{\rm max}$ be an arbitrary element and set 
$a_t\:=\Exp(itX)$ for all $t>0$. For each $t>0$ consider the possibly
unbounded linear functional  

$$f_t\:  {\cal H}_\lambda \to \C, \ \ w\mapsto \int_H \la 
\pi_\lambda(h).v_\lambda, \pi_\lambda(a_t).w\ra \
d\mu_H(h).$$  
In view of Proposition II.4, we have
to show that $\lambda +\rho\in {\cal E}_\Omega$ is equivalent to  
$f_t\in {\cal H}_\lambda'$ for all $t>0$. 

Since $v_\lambda$ is fixed by $H^0$ (cf.\ Proposition II.7(i)), Step 1
and the integration formula of Lemma II.15(ii) yield

$$\eqalign{&\int_H \la \pi_\lambda(h).v_\lambda, \pi_\lambda(a_t).w\ra \ d\mu_H(h)\cr
&\quad =\int_{K^c\cap(HAN)}\la 
\pi_\lambda(h_H(k)).v_\lambda, \pi_\lambda(a_t).w\ra a_H(k)^{-2\rho}\ d\mu_{K^c}(k)\cr
&\quad =\int_{K^c\cap(HAN)}\la 
\pi_\lambda(k a_H(k)^{-1}).v_\lambda, \pi_\lambda(a_t).w\ra a_H(k)^{-2\rho}\ d\mu_{K^c}(k)\cr
&\quad =\int_{K^c\cap(HAN)}\la 
\pi_\lambda(a_t k).v_\lambda, w\ra a_H(k)^{-(\lambda+2\rho)}\ d\mu_{K^c}(k).\cr}\leqno(2.1)$$
Recall from [FH\'O94, Prop.\ 5.3] that 

$${\cal E}_\Omega=\{\lambda\in \a_\C^*\: \int_{K^c\cap(HAN)}
a_H(k)^{-\Re(\lambda+\rho)}\ d\mu_{K^c}(k)
<\infty\}.\leqno(2.2)$$
\par In view of [KN\'O98, Lemma III.15(ii)], the set
$X_t\:=\oline{a_t(K^c\cap HAN)}$ is a compact subset of $HAN$. In
particular we find compact sets $C_H^t$, $C_A^t$, $C_N^t$ contained in
$H$, $A$ and $N$, respectively, such that $X_t\subeq C_H^t C_A^t C_K^t$.
Thus we conclude from Step 1 that 

$$(\forall w\in {\cal H}_\lambda)(\forall x\in X_t)\quad |\la
\pi_\lambda(x).v_\lambda, w\ra|\leq \sup_{a\in C_A^t} a^\lambda
\|v_\lambda\|\cdot\|w\|<\infty.
\leqno(2.3)$$
Hence, in view of (2.1), (2.2) and (2.3) the proof of Step 2 will be
complete, provided we can show that for each element $x$ in the
compact space $X_t$ we can find an  open neighborhood $U\subeq X_t$ of $x$
and an element $w\in {\cal H}_\lambda$ such that $\inf_{y\in U}|\la
\pi_\lambda(y).v_\lambda, w\ra|>0$ holds. But this follows from 
$\la \pi_\lambda(y).v_\lambda, w\ra =\la F(y), w\ra$ and the continuity of
$F$.

\par\nin Step 3: If the integral exists, then its value is $ {\la v_\lambda, v_\lambda\ra
\over \la\nu, v_\lambda\ra} c(\lambda+\rho) \nu$.

\par By Step 1 we know that $\lambda+\rho\in {\cal E}_\Omega$ in the
case where the integral exists. Since $\lambda$ is a highest weight
for an $H\cap K$-spherical representation of $K$, it has to be dominant integral
with respect to $\Delta_k^+$, i.e., $\la \lambda, \alpha\ra \in \N_0$ for
all $\alpha\in \Delta_k^+$. In particular $c(\lambda +\rho)$ exists
(cf.\ Theorem II.14).
Now by Step 2, we know that $\int_H
\pi_\lambda(h).v_\lambda\ d\mu_H(h)\in ({\cal H}_\lambda^{-\omega})^H$. Since 
$\dim ({\cal H}_\lambda^{-\omega})^H\leq 1$ (cf.\ Proposition
II.7(ii)), it follows that $\int_H
\pi_\lambda(h).v_\lambda\ d\mu_H(h)=c\nu$ for some constant $c\in \C$. To
determine $c$ we apply the integral to the element $v_\lambda$. 
With Step 1 and the integration formula of Lemma II.15(i) we 
compute  

$$\eqalign{&\int_H \la\pi_\lambda(h).v_\lambda, v_\lambda\ra 
\ d\mu_H(h)\cr
&\quad =\int_{\oline N\cap (HAN)} \la\pi_\lambda(h_H(\oline
n)).v_\lambda,v_\lambda\ra 
a_H(\oline n)^{-2\rho}\ d\mu_{\oline N} (\oline n)\cr
&\quad =\int_{\oline N\cap (HAN)} \la\pi_\lambda(\oline n a_H(\oline
n)^{-1}).v_\lambda,v_\lambda\ra 
a_H(\oline n)^{-2\rho}\ d\mu_{\oline N} (\oline n)\cr
&\quad = \int_{\oline N\cap (HAN)} \la\pi_\lambda(\oline n)
.v_\lambda, v_\lambda\ra 
a_H(\oline n)^{-(\lambda+2\rho)}\ d\mu_{\oline N} (\oline n)\cr
&\quad =\la v_\lambda, v_\lambda\ra \int_{\oline N\cap (HAN)} 
a_H(\oline n)^{-(\lambda+2\rho)}\ d\mu_{\oline N} (\oline n)\cr
&\quad =\la v_\lambda, v_\lambda\ra c(\lambda+\rho).\cr}$$
This proves Step 3 and completes the proof of the theorem. \qed

\sectionheadline{III. Representations of the relative discrete series}

In this section we state and prove the Harish-Chandra - Godement
Orthogonality relations for homogeneous spaces carrying 
an invariant measure. Then we give the definition of the formal 
dimension $d(\lambda)$ of a unitary highest weight representation
$(\pi_\lambda, {\cal H}_\lambda)$ 
which belongs to the relative discrete series of $H\bs G$. Finally we
derive the formula for $d(\lambda)$ for large values of $\lambda$.

\subheadline{ Orthogonality Relations}

\Definition III.1. Let $G$ be a Lie group, $Z$ its center and $\hat Z$
the group of unitary characters of
$Z$.Let $H\subeq G$ be a closed subgroup. Suppose that $HZ$ is
closed and that $HZ\bs G$ carries an invariant positive measure
$\mu_{HZ\bs G}$. For a fixed $\chi\in \hat Z$ we consider the Hilbert
space of sections
$$\eqalign{\Gamma_\chi^2(H\bs G)&=\{ f\: H\bs G\to\C\: f\
\hbox{measurable}, (\forall z\in Z)(\forall g\in G)\ f(Hzg)=\chi(z)f(Hg);\cr
&\qquad \la f, f\ra_\chi\:=\int_{HZ\bs G} |f(Hg)|^2\ d\mu_{HZ\bs G}(HZg)<\infty\}.\cr}$$
Let $(\pi, {\cal H})$ be an irreducible  unitary 
$H$-spherical representation of $G$ with central character $\chi$. 
Then for all $\nu\in ({\cal H}^{-\omega})^H$  and $v\in {\cal H}^\omega$ we define
a continuous section by 
$$\pi_{v,\nu}\:  H\bs G\to \C, \ \ Hg\mapsto \oline{\la \nu, \pi(g).v\ra}.$$
We say that $(\pi, {\cal H})$ belongs to the {\it relative discrete series
of} $H\bs G$, if there exists non-zero elements $\nu\in ({\cal
H}^{-\omega})^H$ and $v\in {\cal H}^\omega$ such that $\pi_{v,\nu}$
belongs to $\Gamma_\chi^2(H\bs G)$. 
We denote $({\cal H}^{-\omega})_2^H$ the subspace of $({\cal
H}^{-\omega})^H$ which corresponds to the relative discrete series for
$H\bs G$. \qed

In the proof of the following Proposition we adapt a nice idea of
J. Faraut to our setting (cf.\ [Gr96, Sect.\ III.3]). 

\Proposition III.2. {\rm (Orthogonality Relations)} Let $G$ be a Lie
group with center $Z$. Then, if $H$ is a closed subgroup of $G$ such that
$HZ$ is closed and $HZ\bs G$ carries a positive $G$-invariant measure,
then the following assertions hold:

\item{(i)} If $(\pi, {\cal H})$ belongs to the relative discrete
series of $H\bs G$ transforming under the central character $\chi\in
\hat Z$ and $0\neq \nu\in({\cal H}^{-\omega})_2^H$, then all
matrix coefficients $\pi_{v,\nu}$, $v\in {\cal H}^\omega$, belong to 
$\Gamma_\chi^2(H\bs G)$ and there exists a constant $d(\pi,\nu)$ 
depending on the equivalence
class of $\pi$ and on $\nu$ such that the mapping 
$$T\: {\cal H}^\omega\to \Gamma_\chi^2(H\bs G), \ \ v\mapsto 
\sqrt{d(\pi, \nu)} \pi_{v,\nu}$$ 
extends to a $G$-equivariant isometry.

\item{(ii)} If $(\pi, {\cal H})$ and $(\sigma, {\cal K})$ are  two
inequivalent representations of the relative discrete series of $HZ\bs
G$ transforming under the same central character for $Z$, then 
for $\nu\in ({\cal H}^{-\omega})_2^H$ and  $\eta\in ({\cal K}^{-\omega})_2^H$
one has
$$\la \pi_{v,\nu}, \sigma_{w,\eta}\ra=\int_{HZ\bs G} \la \nu, \pi(g).v\ra \oline{\la \eta, \sigma(g).w\ra}
\ d\mu_{HZ\bs G} (HZg) =0$$
for all $v\in {\cal H}^\omega$ and $w\in {\cal K}^\omega$.

\Proof. (i) (cf.\ [Gr96, Sect.\ III.3]) 
Let $D\:=\{ v\in {\cal H}^\omega\: \pi_{v,\nu}\in
\Gamma_\chi^2(H\bs G)\}$ and consider the unbounded operator 
$$S\:D\to \Gamma_\chi^2(H\bs G), \ \ v\mapsto \pi_{v,\nu}.$$
Since $\mu_{HZ\bs G}$ is $G$-invariant by assumption, the same holds
for $D$ and therefore $D$ is dense in ${\cal H}$ by the 
irreducibility of $(\pi,{\cal H})$. We define a positive hermitian 
form on $D$ by 

$$(v|w)\:=\la v, w\ra +\la S.v, S.w\ra_\chi\leqno(3.1)$$
for $v,w\in D$. Denote by $\oline D$ the Hilbert completion of
$D$ with respect to $(\cdot|\cdot)$ and denote the extension of
$(\cdot|\cdot)$
to its completion by the same symbol. Since $\oline D$ is
continuously embedded into ${\cal H}$, there exists a bounded  selfadjoint
injective operator $A\in B({\cal H})$ such that $\im A=\oline D$ and 
$(A.v|w)=\la v,w\ra$ for all $v\in {\cal H}$, $w\in \oline D$.
Since $\la \cdot, \cdot\ra_\chi$ is $G$-invariant by the
$G$-invariance of $\mu_{HZ\bs G}$, it follows from (3.1) that $A$
commutes with $\pi(G)$. Now Schur's Lemma applies and yields 
$A=c \id$ for some constant $c>0$. Thus we deduce from (3.1) that 
$$\la S.v, S.w\ra_\chi=\big({1\over c}-1\big)\la v, w\ra$$
for all $v, w\in D$. In particular $d(\pi, \nu)\:=\big({1\over
c}-1\big)>0$. Moreover $S$ being weakly continuous, its extension to
${\cal H}^\omega$ coincides with ${1\over \sqrt{d(\pi, \nu)}}T$, concluding the proof of (i).
   
\par\nin (ii) Let $T_\pi\: {\cal H}\to  \Gamma_\chi^2(H\bs G)$ and $T_\sigma\: 
{\cal K}\to \Gamma_\chi^2(H\bs G)$ be the equivariant isometric embeddings from (i). 
If $\im T_\pi\cap\im T_\sigma\neq \{0\}$, then 
$$T_\sigma^*\circ T_\pi\: {\cal H}\to {\cal K}$$
describes a non-trivial $G$-equivariant map. By Schur's Lemma $T_\sigma^*\circ T_\pi$
is a scalar multiple of an isometric isomorphism, contradicting 
the inequivalence of $(\pi, {\cal H})$ and $(\sigma, {\cal K})$.\qed

\Remark III.3. If $H\bs G$ is a semisimple symmetric space, then 
the space $({\cal H}^{-\omega})_2^H=({\cal H}^{-\infty})_2$ is finite dimensional 
(cf.\ [Ba87, Th.\ 3.1]). Then Proposition III.2(i)
says that one can find an inner product on $({\cal H}^{-\omega})_2^H$ such that 
$$({\cal H}^{-\omega})_2^H\otimes {\cal H}^\omega\to \Gamma_\chi^2(H\bs G), \ \ 
\nu\otimes v\mapsto \sqrt{d(\pi,\nu)}\pi_{v,\nu}$$
extends to a $G$-equivariant isometry (with $G$ acting trivially on the first 
factor $({\cal H}^{-\omega})_2^H$ of the tensor product). \qed

\subheadline{The formal dimension}

If $G$ denotes a unimodular locally compact group and $L\subeq G$ a closed unimodular
subgroup, then we denote by $\mu_{L\bs G}$ a positive right
$G$-invariant measure on the homogeneous space $L\bs G$.

\Definition III.4. Let $(\pi_\lambda, {\cal H}_\lambda)$ be an
$H$-spherical unitary highest weight representation of $G$ and $0\neq \nu\in
({\cal H}_\lambda^{-\omega})^H$. If $v_\lambda$ is a highest weight vector
for $(\pi_\lambda, {\cal H}_\lambda)$, then the {\it formal
dimension} $d(\lambda)$ of $(\pi_\lambda, {\cal H}_\lambda)$ is defined by 

$${1\over d(\lambda)}\:={1\over |\la \nu, v_\lambda\ra|^2 }\int_{HZ\bs
G} |\la \nu, \pi_\lambda(g).v_\lambda\ra|^2\ d\mu_{HZ\bs G} (HZg).$$
Recall that $\la \nu, v_\lambda\ra \neq 0$ and that the definition of $d(\lambda)$ is independent of the
particular choice of $v_\lambda$ and $0\neq \nu\in ({\cal
H}_\lambda^{-\omega})^H$ (cf.\ Proposition II.7(ii)). 
\par The relation between the number $d(\pi_\lambda,\nu)$ from Proposition
III.2 and $d(\lambda)$ is given by $d(\lambda)={|\la \nu,
v_\lambda\ra|^2\over \la v_\lambda, v_\lambda\ra } d(\pi_\lambda,
\nu)$. In particular, if $\nu$ is normalized by ${|\la \nu,
v_\lambda\ra|^2\over \la v_\lambda, v_\lambda\ra}=1$, then we have
$d(\lambda)=d(\pi_\lambda, \nu)$. \qed

\Remark III.5. The particular normalization of $d(\lambda)$ as in 
Definition III.4 is motivated from Harish-Chandra's treatment of the ``group case''
(cf.\ [HC56]). The group case is defined by 
$G=G_0\times G_0$ and $H=\Delta(G)=\{(g,g)\: g\in G_0\}$ for a simply connected hermitian 
Lie group $G_0$. 
Then we have a natural isomorphism 
$$G_0\to H\bs G, \ \ g\mapsto H(g,1)$$
and the invariant measure $\mu_{ZH\bs G}$ corresponds to a
Haar measure  $\mu_{Z(G_0)\bs G_0}$ on $Z(G_0)\bs G_0$.

\par The spherical unitary highest weight representations of $G$ are 
given by $(\pi_\lambda\otimes\pi_\lambda^*, {\cal H}_\lambda\hat\otimes 
{\cal H}_\lambda^*)$ with $(\pi_\lambda, {\cal H}_\lambda)$ a unitary highest weight representation 
of $G_0$ and $(\pi_\lambda^*, {\cal H}_\lambda^*)$ its dual representation. 
Recall that $ {\cal H}_\lambda\hat\otimes {\cal H}_\lambda^*$ is isomorphic to 
the space of Hilbert-Schmidt operators $B_2({\cal H}_\lambda)$ on 
${\cal H}_\lambda$ and that the corresponding analytic vectors are of trace class, i.e., 
$B_2({\cal H}_\lambda)^\omega\subeq B_1({\cal H}_\lambda)$ (cf.\ [HiKr99a, App.]). The up to 
scalar unique $H$-fixed hyperfunction vector is given by the conjugate trace: 
$$\nu\: B_2({\cal H}_\lambda)^\omega \to \C, \ \ A\mapsto \oline {\tr (A)}.$$  
Further a highest weight vector for $(\pi_\lambda\otimes\pi_\lambda^*, {\cal H}_\lambda\hat\otimes 
{\cal H}_\lambda^*)$ is given by $v_\lambda\otimes v_\lambda^*$. Then 
$\la \nu, v_\lambda\otimes v_\lambda^*\ra =\la v_\lambda, v_\lambda\ra $ and the expression 
for $d(\lambda)$ from Definition III.4 gives that 
$${1\over d(\lambda)}={1\over |\la v_\lambda, v_\lambda\ra|^2 }\int_{Z(G_0)\bs G_0}
|\la \pi_\lambda(g).v_\lambda, v_\lambda\ra|^2\ d\mu_{Z(G_0)\bs G_0} (Zg).$$
Thus we see that our definition of the formal dimension coincides in the group 
case with the standard one introduced by  Harish-Chandra (cf.\ [HC56]). \qed

\Theorem III.6. Let $(\pi_\lambda, {\cal H}_\lambda)$
be an unitary highest weight representations of
$G$ for which $(\pi_\lambda^K, F(\lambda))$ is $H\cap K$-spherical. 
Assume that $\lambda+\rho\in {\cal
E}_\Omega$ and that $(\pi_\lambda, {\cal H}_\lambda)$ belongs to the
holomorphic discrete series of $G$. Then $(\pi_\lambda, {\cal H}_\lambda)$ is
$H$-spherical, belongs to the relative discrete series of $H\bs G$ and
the formal degree $d(\lambda)$ is given by
$$d(\lambda)=d(\lambda)^G c(\lambda+\rho),$$
where $d(\lambda)^G$ is the formal dimension of $(\pi_\lambda, {\cal H}_\lambda)$ relative to $G$.

\Proof. Since $\lambda+\rho\in {\cal E}_\Omega$ the asumptions of 
Theorem II.16 are satisfied and the theorem applies. 
Thus $(\pi_\lambda, {\cal H}_\lambda)$ is
$H$-spherical and if $0\neq \nu\in ({\cal
H}_\lambda^{-\omega})^H$  and $v_\lambda$ is a
highest weight vector, then we have 

$$\nu={\la \nu, v_\lambda\ra  \over \la v_\lambda, v_\lambda\ra c(\lambda +\rho)} 
\int_H \pi_\lambda(h).v_\lambda \ d\mu_H(h).\leqno(3.2)$$
If we insert (3.2) in the formula for $\nu$ in the definition of the formal
dimension we obtain that 

$$\eqalign{{1\over d(\lambda)}&={1\over |\la \nu, v_\lambda\ra|^2}\int_{HZ\bs G} 
|\la \nu, \pi_\lambda(g).v_\lambda\ra|^2\ d\mu_{HZ\bs G} (HZg)\cr
&={1\over \la v_\lambda, v_\lambda\ra^2 c(\lambda+\rho)^2} 
\int_{HZ\bs G}\int_H\int_H \la
\pi_\lambda(h_1).v_\lambda, 
\pi_\lambda(g).v_\lambda\ra\cr 
&\qquad\qquad \la \pi_\lambda(g).v_\lambda,
\pi_\lambda(h_2).v_\lambda\ra
\ d\mu_H(h_1)\ d\mu_H(h_2)\ d\mu_{HZ\bs G} (HZg)\cr
&={1\over\la v_\lambda, v_\lambda\ra^2 c(\lambda+\rho)^2} 
\int_{HZ\bs G}\int_H\int_H \la
\pi_\lambda(h_2h_1).v_\lambda, 
\pi_\lambda(g).v_\lambda\ra\cr 
&\qquad\qquad \la \pi_\lambda(h_2^{-1}g).v_\lambda, v_\lambda\ra
\ d\mu_H(h_1)\ d\mu_H(h_2)\ d\mu_{HZ\bs G} (HZg)\cr
&={1\over \la v_\lambda, v_\lambda\ra^2 c(\lambda+\rho)^2} 
\int_{HZ\bs G}\int_H\int_H \la
\pi_\lambda(h_1).v_\lambda, 
\pi_\lambda(h_2g).v_\lambda\ra\cr 
&\qquad\qquad \la \pi_\lambda(h_2g).v_\lambda,
v_\lambda\ra
\ d\mu_H(h_1)\ d\mu_H(h_2)\ d\mu_{HZ\bs G} (HZg)\cr
&={1\over \la v_\lambda, v_\lambda\ra^2 c(\lambda+\rho)^2} 
\int_H\int_{HZ\bs G}\int_H \la
\pi_\lambda(h_1).v_\lambda, 
\pi_\lambda(h_2g).v_\lambda\ra\cr 
&\qquad\qquad \la \pi_\lambda(h_2g).v_\lambda,
v_\lambda\ra
\ d\mu_H(h_2)\ d\mu_{HZ\bs G} (HZg)\ d\mu_H(h_1)\cr
&={1\over \la v_\lambda, v_\lambda\ra^2 c(\lambda+\rho)^2} 
\int_H\int_{Z\bs G}\la
\pi_\lambda(h_1).v_\lambda, 
\pi_\lambda(g).v_\lambda\ra\
\la \pi_\lambda(g).v_\lambda,
v_\lambda\ra
\ d\mu_{Z\bs G} (Zg)\ d\mu_H(h_1).\cr}$$
Thus if we apply the Harish-Chandra-Godement Orthogonality Relations
for $L^2(Z\bs G)$ and once more (3.2) we obtain that  

$$\eqalign{{1\over d(\lambda)}&={1\over d(\lambda)^G}\cdot{1\over \la v_\lambda,
v_\lambda\ra^2 c(\lambda+\rho)^2}
\la v_\lambda, v_\lambda\ra \int_H \la \pi_\lambda(h).v_\lambda,
v_\lambda\ra
\ d\mu_H(h)\cr
& ={1\over d(\lambda)^G}\cdot{1\over \la v_\lambda, v_\lambda\ra
c(\lambda+\rho)^2} c(\lambda+\rho) \la v_\lambda, v_\lambda\ra  
={1\over d(\lambda)^G c(\lambda+\rho)},\cr}$$  
as was to be shown.\qed

\sectionheadline{IV. Analytic continuation in $\lambda$}

In this section we prove the analytic continuation of the formula
for the formal dimension $d(\lambda)$ from Theorem III.6. The proof 
is quite technical and we need some preparation on algebraic and
analytic level.

\subheadline{Algebraic preliminaries}

In this subsection we collect some facts concerning the 
fine structure theory of compactly 
causal symmetric Lie algebras. The results are mainly
due to \'Olafsson (cf.\ [\'Ol91]).

\Lemma IV.1. Let $(\g,\tau)$ be a compactly causal symmetric Lie
algebra, then we can choose root vectors $E_\alpha \in \g_\C^\alpha$,
$\alpha\in\hat\Delta_n$, such that the following conditions are
satisfied:
\item{(1)} $\oline{E_\alpha }= E_{-\alpha}$.
\item{(2)} $\alpha(H_\alpha)=2$ with $H_\alpha=[E_\alpha ,
E_{-\alpha}]$.
\item{(3)} $\tau (E_\alpha )= E_{\tau\alpha}$, where $\tau\alpha=\tau\circ\alpha$.

\Proof. Let $\kappa$ denote the Cartan-Killing form on $\g_\C$ and
define a hermitian inner product on $\g_\C$ by $\la X, Y\ra\:=-\kappa(X,
\theta(\oline Y))$.

\par For each $\alpha\in\hat\Delta_n^+$ let $0\neq E_\alpha \in
\g_\C^\alpha$ be an arbitrary element of length $1$. Then define $E_{-\alpha}$ by
$E_{-\alpha}\:=\oline{E_\alpha }$. Thus (1) is satisfied. 
Now $\tau(E_\alpha )\subeq \C E_{\tau\alpha}$ implies the existence of
complex numbers $c_\alpha$ such that $\tau(E_\alpha )=c_\alpha
E_{\tau\alpha}$. Now $\tau$ being an involutive implies $c_\alpha
c_{\tau\alpha}=1$, further $\tau$ being an isometry implies that 
$|c_\alpha|=1$ and finally $\tau$ being complex linear implies that
$\oline{c_\alpha}=c_{-\alpha}$ for all $\alpha\in \hat\Delta_n$.  
Thus $c_{\tau\alpha}= \oline{c_\alpha}=c_{-\alpha}$. 
For each complex number $z=e^{i\phi}$, $\phi\in [0, 2\pi[$, of
modulus $1$ we define $z^{1\over 2}=e^{i{\phi\over 2}}$. 
Thus redifining $E_\alpha $, $\alpha\in \hat\Delta_n^+$,  by $\oline {c_\alpha}^{1\over 2}
E_\alpha $, leaves (1) untouched and in addition satisfies (3).
\par Since $\g_\C^\alpha\subeq \p_\C$ for all $\alpha\in \hat\Delta_n^+$, we have
$\alpha([E_\alpha , E_{-\alpha}])>0$, and so by rescaling $E_\alpha $ with
an appropriate positive number we may in addition assume that (2)
holds. This proves the lemma.\qed

Let $\hat\Gamma=\{\hat\gamma_1, \ldots, \hat\gamma_r\}$ be a maximal system of
strongly orthogonal, i.e., $\hat\gamma_j \pm \hat\gamma_i$ is never a
root and $\hat\Gamma\subeq \hat\Delta_n^+$ has maximal many elements  with respect to this
property. In view of [Hi\'Ol96, Lemma 4.1.7] or [\'Ol91, Sect.\ 3], we may choose $\hat\Gamma$
invariant under $-\tau$.

For each $1\leq j\leq r$ we set $\hat E_j\:=E_{\hat \gamma_j}$, 
$\hat E_{-j}\:=E_{-\hat\gamma_j}$ and  
$\hat X_j\:=i(\hat E_j -\hat E_{-j})$. According to [HC56, Cor.\ to Lemma 8], the space 

$$\e\:=\bigoplus_{j=1}^r \R \hat X_j=\bigoplus_{j=1}^r \R i(\hat E_j -
\hat E_{-j})$$
is maximal abelian in $\p$. Note that $\e$ is $\tau$-invariant by
the special choice of the non-compact root vectors  (cf.\ Lemma IV.5(3)) and the
$-\tau$-invariance of $\hat\Gamma$.

We consider the {\it Cayley transform}
$$C=e^{i{\pi\over 4}\ad (\sum_{j=1}^r \hat E_j + \hat E_{-j})}$$
which is an automorphism of $\g_\C$. Finally we set $\hat H_j\:=
H_{\hat\gamma_j}$ for all $1\leq j\leq r$.

\Lemma IV.2. The Cayley transform $C$ has the following properties:

\item{(i)} For all $1\leq j\leq r$ one has $C(\hat X_j)= \hat H_j$ and
$C(\hat H_j)=-\hat X_j$.
\item {(ii)} We have $i{\pi\over 4}(\sum_{j=1}^r \hat E_j +
\hat E_{-j})\in  i \h_\p$. In particular, one has 
\itemitem{(a)} $\tau\circ C= C\circ \tau$,
\itemitem{(b)} $\theta\circ C= C^{-1}\circ \theta$,
\item{(iii)} The Cayley transform yields an isomorphism $C\:
\e\to C(\e)$ with $C(\e)\subeq i \t$ a $\tau$-invariant subspace. 

\Proof. (i) This follows from $\sL(2,\R)$-reduction (cf.\ [HC56, p.\
584], [Hi\'Ol96, Lemma A.3.2(3)]).
\par\nin (ii) It follows from $\g_\C^{\hat\alpha}\subeq \p_\C$, for all
$\hat\alpha\in \hat\Delta_n$ and  Lemma IV.1(1) that $i{\pi\over
4}(\sum_{j=1}^r \hat E_j +
\hat E_{-j})\in i\p$. Further Lemma IV.1(3) and the 
$-\tau$-invariance of $\hat\Gamma$ imply
$$\tau\big(\sum_{j=1}^r \hat E_j+\hat E_{-j}\big)=\sum_{j=1}^r
\tau(\hat E_j)+ \tau(\hat E_{-j})=\sum_{j=1}^r E_{\tau\hat\gamma_j}+ 
E_{-\tau\hat\gamma_j}=\sum_{j=1}^r \hat E_j+ \hat E_{-j}.$$
Thus $i{\pi\over 4}(\sum_{j=1}^r \hat E_j +
\hat E_{-j})\in i\h_\p$. This proves (i).
\par\nin (iii) This follows from (i) and (ii)(a).\qed

Recall that $\e$ is $\tau$-invariant and write $\b=\e\cap\q$ for the
set of $-\tau$-fixed points.

\Lemma IV.3. Let $\cc\:=C(\b)$. Then $\cc\subeq \a$
and the Cayley transform yields an isomorphism $C\: \b\to\cc$. 

\Proof. Since $C(\b)\subeq i\t$ by Lemma IV.2(i), the fact that
$\b\subeq \q$ and that $C$ commutes with $\tau$ (cf.\ Lemma IV.2(ii))
imply that $C(\b)\subeq i(\t\cap\q)$. But $i(\t\cap\q)=\a$ by the
definition of $\a$, proving the lemma.\qed

Recall that $\b$ is maximal abelian
subspace of $\q\cap\p$ (this follows for instance from the $c$-dual 
version of Lemma 4.1.9 in [Hi\'Ol96]) and denote by $\Sigma=\Sigma(\g,
\b)$ the set of roots of $\g$ with respect to $\b$. Recall that 
$\Sigma$ is an abstract root system (cf.\ [Sch84, Sect.\ 7.2]). We write 

$$\g=\z_\g(\b)\oplus\bigoplus_{\phi\in \Sigma}\g^\phi$$
for the corresponding root space decomposition. 
By  Lemma IV.3, 
the Cayley transform induces a mapping 
$C^t\: \a^*\to \b^*, \ \ \alpha\mapsto \alpha\circ C\res_\b$ 
and we set 
$$\Sigma_n=C^t(\Delta_n)\res_\b\quad \hbox{and}\quad \Sigma_k=
C^t(\Delta_k)\res_\b\bs\{0\}. $$
Let $\Gamma=\{{1\over 2}(\hat\gamma_j -\tau\hat\gamma_j)\: 1\leq j\leq r\}$ denote the
restricted set of strongly orthogonal roots. Note that $\Gamma\subeq
\cc^*$ by Lemma IV.2(i). Thus we can write 
$\Gamma=\{\gamma_1, \ldots, \gamma_s\}$ for some $1\leq s\leq r$. For
each $1\leq j\leq s$ we define $H_j\in \cc$ by $\gamma_j(H_j)=2$ and 
$\gamma_k(H_j)=0$ for $k\neq j$. We set $X_j\:=-C(H_j)$ 
for all $1\leq j\leq s$. Then  
$$\b=\bigoplus_{j=1}^s \R X_j.$$
As a final algebraic tool we need explicit information
on the root system $\Sigma$ which is provided  by \'Olafsson's Theorem on 
double restricted root systems (cf.\ [\'Ol91, Sect.\
3], [H\'O\O91, Prop.\ 3.1]). 
For all $1\leq j\leq s$ we set $\psi_j\:=C^t(\gamma_j)$ and
note that $\psi_j(X_j)=2$ since $C(X_j)=H_j$ (cf.\ Lemma IV.2(i), (ii)).  

Finally we put $\Sigma^+\:=C^t(\Delta^+)\res_\b\bs \{0\}$, $\Sigma_n^+\:=\Sigma_n\cap\Sigma^+$
and $\Sigma_k^+\:=\Sigma_k\cap\Sigma^+$.

\Theorem IV.4. {\rm (\'Olafsson)} If $(\g,\tau)$ is compactly causal, 
then the following
assertions concerning the double restricted root system $\Sigma=\Sigma(\g,\b)$ hold:
\item{(i)} The restricted root system has the following form 
$$\Sigma_k=\pm\{ {1\over 2}(\psi_i
-\psi_j): i< j\}\cup\pm\{{1\over 2}\psi_j\: 1\leq j\leq s\} $$
and
$$\Sigma_n^+=\{ {1\over 2}(\psi_i +\psi_j)\: 1\leq i, j\leq
s\}\cup\{{1\over 2}\psi_j\: 1\leq j\leq s\}. $$
The second sets in $\Sigma_k$ and $\Sigma_n^+$ are empty if and only if $C^4=\id$. 
If further $\psi_s$ is chosen to be a simple root, then 
$$\Sigma_k^+\subeq \{ {1\over 2}(\psi_i
-\psi_j): i< j\}\cup\{{1\over 2}\psi_j\: 1\leq j\leq s\}.$$
\item{(ii)} All $\psi_j$, $1\leq j\leq s$, have the same length.  
\item{(iii)} The conjugacy classes of the restricted root system
under the Weyl group associated to $\Sigma$ are given by 
\itemitem{(1)} $\{\pm{1\over 2}(\psi_i\pm\psi_j)\: 1\leq i, j\leq
s, i\neq j\}$ 
\itemitem{(2)} $\{\pm\psi_j\: 1\leq  j\leq s\}$ 
\itemitem{(3)} $\{\pm{1\over 2}\psi_j\: 1\leq j\leq
s\}$

\Proof. (i) Let $\hat\Sigma=\hat\Sigma(\g,\e)$ be the restricted root system with
respect to the maximal abelian subspace $\e$ and $\hat\Sigma_k$,
$\hat\Sigma_n$ defined as above. Write
$\hat\psi_j\:=C^t(\hat\gamma_j)$ for all $1\leq j\leq r$.
Suppose first that $\g$ is simple. Then for the analogous statement 
for $\hat\Sigma$ in stead of $\Sigma$ and $\hat\psi_j$ in stead of 
$\psi_j$, Harish-Chandra has proved in [HC56, Lemma 13-16] that 
$\hat \Sigma_k,\hat\Sigma_n^+$  are contained in the
asserted subsets, Moore proved equality (cf.\ [Mo64, Th.\ 2]) and 
finally Kor\'anyi and Wolf have shown in [KoWo65, Prop.\ 4.4 with Remark]
that the second set in $\hat\Sigma_n^+$ is empty if and only if
$C^4=\id$. Now taking restrictions to $\cc$  yields (i)  for $\g$
simple. 
\par In the group case similar considerations lead to the same
result. 
\par\nin (ii) This can be deduced from [Mo64, Th.\ 2(2)], but we 
propose here a much simpler proof. We use (i) and the fact that $\Sigma$ is an abstract
root system. As usual we write $s_\psi$, $\psi\in
\Sigma$, for the reflection associated to $\psi$. Then we obtain 
for all $1\leq i\not= j\leq s$ that 
$$\eqalign{s_{{1\over2}(\psi_i+\psi_j)}(\psi_j)&= \psi_j 
-{2\la \psi_j, {1\over2}(\psi_i+\psi_j)\ra\over 
\la{1\over2}(\psi_i+\psi_j), {1\over2}(\psi_i+\psi_j)\ra}
{1\over2}(\psi_i+\psi_j)\cr
&=\psi_j -{2\la \psi_j, \psi_j\ra\over 
\la\psi_i, \psi_i\ra +\la \psi_j,\psi_j\ra}(\psi_i+\psi_j).\cr}$$
Thus it follows from (i) and $s_{{1\over2}(\psi_i+\psi_j)}(\psi_j)\in \Sigma$ that 
${\la \psi_j, \psi_j\ra\over 
\la\psi_i, \psi_i\ra +\la \psi_j,\psi_j\ra}\in\{{1\over 2},{1\over
4}\}$. Interchanging $i$ and $j$ then yields ${\la \psi_j, \psi_j\ra\over 
\la\psi_i, \psi_i\ra +\la \psi_j,\psi_j\ra}={1\over 2}$
or equivalently that $\la \psi_j,\psi_j\ra =\la \psi_i, \psi_i\ra$. This
proves (ii). 

\par\nin (iii) In view of (i), we have for all $1\leq i,j,k\leq r$ that
$$\eqalign{& s_{{1\over2}(\psi_i\pm \psi_j)}(\psi_j)=\mp \psi_i,\cr
& s_{{1\over2}(\psi_i\pm\psi_j)}\big({1\over 2}(\psi_j
\pm\psi_k)\big)={1\over 2}(\mp\psi_i\pm\psi_k),\cr
& s_{\psi_i}\big({1\over 2}(\psi_i \pm\psi_j)\big)={1\over 2}(
-\psi_i \pm \psi_j).\cr}\leqno(4.1)$$
This proves (iii).\qed 

From now on we assume that $\psi_s$ is a simple root. 
Then Theorem IV.4(i) says that 
$$\Sigma_n^+=\{ {1\over 2}(\psi_i +\psi_j)\: 1\leq i, j\leq
s\}\cup \{{1\over 2}\psi_j\: 1\leq j\leq s\}.\leqno(4.2)$$

and 
$$\Sigma_k^+=\{{1\over 2}(\psi_i -\psi_j)\: 1\leq i<j\leq s \}\cup
\{{1\over 2}\psi_j\: 1\leq j\leq s\}.\leqno(4.3)$$
Further it follows from Theorem IV.4(i) and the first formula in (4.1)
that the Weyl group ${\cal W}(\Sigma_k)$  of $\Sigma_k$ acts on $\b$
as the full permutation group of
the $X_j$'s. 

We write $\b^+=\{X\in \b\: (\forall \phi\in\Sigma^+)\
\phi(X)\geq 0\}$ for the Weyl chamber corresponding to $\Sigma^+$. By
(4.2) and (4.3) we then have 
$$\b^+=\{\sum_{j=1}^s x_j X_j\: 0\leq x_s\leq \ldots\leq x_1\}.$$
Further let  $\a^+\:=\{ X\in \a: (\forall \alpha\in
\Delta^+)\ \alpha(X)\geq 0\}$ and $\cc^+\:=\a^+\cap \cc$. Note that $C(\b^+)=\cc^+$ by
the construction of $\Sigma^+$.

\Lemma IV.5. The following equality holds
$$C_{\rm min}^\star \cap(-\check C_k^\star)= (\cc^+)^\star\cap(-\check
C_k^\star), $$ 
where the stars $\star$ are all taken in $\a^*$. 

\Proof. First recall some basic rules in dealing with convex cones
(cf.\ [Ne99b, Ch.\ V]). If
$W$ is a closed convex cone in an euclidean space $V$, then
$(W^\star)^\star=W$. Further for two closed convex cones $W_1,
W_2\subeq V$ we have $(W_1\cap W_2)^\star=\oline{W_1^\star+W_2^\star}$.

\par Let now the convex cone on the left hand side be denoted by $W_1$, the
other one by $W_2$. Let $p\: \a\to\cc$ be the orthogonal projection
with respect to the Cartan-Killing form. We claim that
$p(W_1^\star)=p(W_2^\star)$. Assume first that no half roots in 
$\Sigma$ occur. Then from the
Cayley-transform analogs of (4.2) and (4.3) it follows that 

$$p(W_1^\star)=p(\oline{C_{\rm min} -\check C_k})=\bigoplus_{j=1}^s\R^+ H_j+
\bigoplus_{j=1}^{s-1} \R^+(H_{j+1}-H_j),$$
and  
$$p(W_2^\star)=p(\oline{\cc^+ -\check C_k})=\Big(\big(\bigoplus_{j=1}^s\R^+
H_j\big)\cap
\{\sum_{j=1}^s x_j H_j\: x_s\leq \ldots\leq x_1\}\Big)+
\bigoplus_{j=1}^{s-1} \R^+(H_{j+1}-H_j).$$
From these two equalities the claim follows in the case of no half roots in $\Sigma$.
The general case  is easily deduced from this. 

\par Let $r\: \a^*\to\cc^*, \  r(\lambda)\:=\lambda\res_\cc$ be the
restriction map and note that $r$ is the dual map of the inclusion mapping
$\cc\to\a$. Since both $W_1$ and $W_2$
are closed, we have $(W_{1,2}^\star)^\star= W_{1,2}$, and so 
$$W_{1,2}\res_\cc=r(W_{1,2})=(p(W_{1,2}^\star))^\star.$$
Hence our claim
implies that $W_1\res_\cc=W_2\res_\cc$. Thus $W_1\subeq W_2$ by the 
definition  of $W_1$ and $W_2$. 
\par For the converse inclusion we first note that an element $\lambda\in 
-\check C_k^\star$ belongs to $W_1$ if and only if
$\lambda(\check\beta)\geq 0$, where $\beta$ is the highest root (this
becomes clear from our construction of the positive systems). Recall
that $\hat\Gamma$ can be constructed inductively starting with the
highest root (cf.\ [HC56, p.\ 108]). Thus 
$\beta=\gamma_1in \Gamma$. Hence  $W_1=(\gamma_1)^\star \cap -\check
C_k^\star$, and so $W_2\subeq W_1$ since $(\cc^+)^\star \subeq
(\gamma_1)^\star$. \qed

\subheadline{Analytic preliminaries}

Recall the definition of $\b^+$ and set $B^+\:=\exp(\b^+)$.

\Lemma IV.6. {\rm (Flensted-Jensen)} Let $L=Z_{H\cap K}(\b)$. Then for the homegeneous space
$HZ\bs G$ the following assertions hold:

\item{(i)} The subgroups $HZ$ and $L Z$ of $G$ are closed and $Z\bs
L Z$ is compact. 
\item{(ii)} The mapping 
$$\Phi\: B^+\times L Z\bs K\to HZ\bs G,\ \  (b, LZk)\mapsto LZbk$$
is a diffeomorphism onto its open image. The image is dense with
complement of Haar measure zero. 
\item{(iii)} Up to normalization of measures we have for all 
$f\in L^1(HZ\bs G)$ the following integration formula
$$\int_{HZ\bs G} f(HZg) \ d\mu_{HZ\bs G}(HZg)=
\int_{Z\bs K}\int_{\b^+}f(HZ\exp(X)k)\ J(X)\ dX\ d\mu_{Z\bs K}(Zk)$$
with 
$$J(X)=\prod_{\phi\in\Sigma^+}\cosh(\phi(X))^{m_\phi^+}\sinh(\phi(X))^{m_\phi^-},$$
where $m_\phi^\pm\:=\dim(\{ X\in \g^\phi\: \theta\tau(X)=\pm X\})$.

\Proof. (i) The closedness of $HZ$ and $L Z$ follows from the
closedness of $\Ad(H)$ and $Z_{\Ad(H)}(\b)$ in the adjoint group
$\Ad(G)$. Finally $Z\bs L Z$ is a closed subgroup of the compact group
$Z\bs Z(H\cap K)$ and hence compact. 
\par\nin (ii) [Sch84, Prop.\ 7.1.3]. 
\par\nin (iii) It follows from [FJ80, Th.\ 2.6] or [Sch84, Lemma
8.1.2] that 
 
$$J(X)\:=\det\big(d\Phi(X,LZk)\big)=\prod_{\phi\in\Sigma^+}
\cosh(\phi(X))^{m_\phi^+}\sinh(\phi(X))^{m_\phi^-}$$
for all $X\in \b^+$ and $k\in K$. Thus it follows from (ii)
that 

$$\int_{HZ\bs G} f(HZg) \ d\mu_{HZ\bs G}(HZg)=
\int_{LZ\bs K}\int_{\b^+}f(HZ\exp(X)k)\ J(X)\ dX\ d\mu_{LZ\bs K}(Zk)$$
holds for all $f\in L^1(HZ\bs G)$. In view of (i), we may replace the 
integration over $LZ\bs K$ by an  $Z\bs K$-integral, proving
(iii).\qed

\Lemma IV.7. Realize $G$ as a submanifold of $\tilde M\times P^+$ as
in {\rm Proposition II.10(ii)}. Then  for $b=\exp_G(\sum_{j=1}^s x_j X_j)\in B$
and 
$$\mu(b)\:=\exp_{\tilde{K_\C}}\Big(\sum_{j=1}^s {1\over 2} \log\cosh(2x_j)
H_j\Big)\in A\subeq \tilde{K_\C}$$
the following assertions hold: 
\item{(i)} We have $b\in \{ [h, \mu(b)]\: h\in \tilde{H_\C}\} \times P^+$. 
\item{(ii)} If $X\in \b^+$, then $\log\mu(\exp_G(X))\in \cc^+$. 

\Proof. (i) This follows directly from [Hi\'Ol96, p.\ 210-211].   
\par\nin (ii) Recall that $X=\sum_{j=1}^s x_j X_j\in \b^+$ if and only if 
$0\leq x_s\leq \ldots \leq x_1$. Now the assertion follows from (i) and the 
monotonicity of the mapping $\R^+\to \R, \ x\mapsto \log\cosh(x)$.\qed

\subheadline{Proof of the analytic continuation}

Let $(\pi_\lambda, {\cal H}_\lambda)$ be an
$H$-spherical unitary highest weight representation of $G$. Further
let $\nu\in
({\cal H}_\lambda^{-\omega})^H$ an $H$-fixed hyperfunction vector
and $\nu_0=\nu\res_{F(\lambda)}\in F(\lambda)^{H\cap K}$. 
We normalize $\nu$ by setting $\|\nu_0\|=1$ and then $v_\lambda$ by 
$|\la \nu, v_\lambda\ra|=1$. 
Then we have 

$$d(\lambda)=I(\lambda)^{-1}\quad \hbox{with}\quad
I(\lambda)\:=\int_{HZ\bs G} |\la \nu, \pi_\lambda(g).v_\lambda\ra|^2\
d\mu_{HZ\bs G} (HZg).$$

\Definition IV.8. On the non-compactly Riemannian symmetric space 
$K(0)\bs G(0)$ we define the 
{\it spherical function with parameter} $\lambda\in \a_\C^*$  by 

$$\phi_\lambda^0(g)=\int_{K(0)} a_H(gk)^{\lambda-\rho_k} \
d\mu_{K(0)}(k), $$
for all $g\in G(0)$.\qed

\Remark IV.9. Note that if $\lambda\in \a^*$ is the highest weight of an $H\cap
K$-spherical representation $(\pi_\lambda^K, F(\lambda))$ of
$\tilde{K_\C}$, then $\phi_{\lambda+\rho_k}^0$ extends to a holomorphic 
function on $\tilde{K_\C}$ and we have 
$$\phi_{\lambda+\rho_k}^0(k)=\la \pi_\lambda^K(k).\nu_0, \nu_0\ra\leqno(4.4)$$ 
for all $k\in \tilde{K_\C}$ (cf.\ [Hel84, Ch.\ V, Th.\ 4.3]). \qed

\Proposition IV.10. With the notation of {\rm  Lemma IV.7} we have 
$$I(\lambda)={1\over\dim F(\lambda)}\int_{\b^+}
\phi_{\lambda+\rho_k}^0(\mu(\exp_G(X))^2)
\ J(X) \ dX$$
where $J(X)$ is given as in {\rm Proposition IV.6(iii)}.

\Proof. (cf.\ [HC56, p.\ 599], [Gr96, Prop.\ 10])
In the sequel we identify $\b$ with $B$ via the exponential mapping
and for $b=\exp_G(X)\in B^+$ we set $J(b)\:=J(X)$. Then by Lemma IV.6(iii) we have 
$$\eqalign{I(\lambda)&=\int_{HZ\bs G} |\la \nu, \pi_\lambda(g).v_\lambda\ra|^2\
d\mu_{HZ\bs G} (HZg)\cr
&=\int_{Z\bs K}\int_{B^+} |\la \nu, \pi_\lambda(bk).v_\lambda\ra|^2\ J(b)
\ d\mu_B(b)
\ d\mu_{Z\bs K}(k). \cr}\leqno(4.5)$$
In view of Lemma II.10(ii), we can write each element in $b\in B^+$ as 
$([h_\C(b), \mu(b)], p_+(b))\in \tilde M\times P^+$ with
$\mu(b)\in \tilde {K_\C}$. Now the
same consideration as in the proof of Step 1 of Theorem II.16 
yields for all $b\in B^+$ and $k\in K$ that 
$$\eqalign{\la \nu, \pi_\lambda(bk).v_\lambda\ra &=
\la \nu, \pi_\lambda\big(([h_\C(b), \mu(b)],
p_+(b))k\big).v_\lambda\ra\cr
&=\la \nu, \pi_\lambda([\1, \mu(b)k],k^{-1}
p_+(b)k).v_\lambda\ra=\la \nu, \pi_\lambda(\mu(b)k).v_\lambda\ra\cr
&=\la \nu_0, \pi_\lambda^K(\mu(b)k).v_\lambda\ra.\cr}$$
If we insert this expression for the matrix coefficient in (4.5),  use
Schur's Orthogonality Relations for $(\pi_\lambda^K, F(\lambda))$
and the relation $\pi_\lambda^K(\mu(b))^*=\pi_\lambda^K(\mu(b))$ (cf.\ Lemma IV.7), 
we arrive at 

$$\eqalign{I(\lambda)&=\int_{B^+} 
\int_{Z\bs K} |\la \nu_0, \pi_\lambda^K(\mu(b)k).v_\lambda\ra|^2\ J(b)\ d\mu_{Z\bs K}(k)
\ d\mu_B(b)\cr
&={1\over \dim F(\lambda)} \int_{B^+} \la \pi_\lambda^K(\mu(b)).\nu_0, 
\pi_\lambda^K(\mu(b)).\nu_0\ra\ J(b)
\ d\mu_B(b)\cr
&= {1\over \dim F(\lambda)} \int_{B^+} \la \pi_\lambda^K(\mu(b)^2).\nu_0, 
\nu_0\ra\ J(b)\ d\mu_B(b).\cr}$$
Now the assertion of the proposition follows from (4.4).\qed

\Lemma IV.11. Let $V$ be a finite dimensional real vector space,
$W\subeq V$ be an open convex cone,  $\alpha_1,\ldots,\alpha_n, \beta_1, \ldots, \beta_m\in
W^\star\bs\{0\}$ 
and $p_1,\ldots, p_n, q_1\ldots, q_m\in \N_0$. For 
every $\lambda\in V^*$ we define the integral

$$H(\lambda)\:=\int_W e^{\lambda(x)} 
\prod_{j=1}^n (\sinh\alpha_j(x))^{p_j}\prod_{\j=1}^m \big(\cosh\beta_j(x)\big)^
{q_j}\ d\mu_V(x).$$
Then $H(\lambda)$ converges if and only if $\lambda+\sum_{j=1}^np_j\alpha_j
+\sum_{j=1}^m q_j \beta_j\in 
-\Int W^\star$. 

\Proof. If $q_1=\ldots=q_m=0$, then this is Lemma IV.6 in [Kr98]. The
general case is easily obtained from this. \qed

The following characterization of the relative discrete
series by the parameter $\lambda$ is due to Hilgert,
\'Olafsson and \O rsted and was obtained in two steps (cf.\ [\'O\O 91, Th.\
5.2], [H\'O\O 91, Th.\ 3.3]). We present an essentially  modified proof here, but we
point out that it is not our objective to give new proofs of
well-known facts. In the course of our arguments we obtain an important new estimate which is 
crucial for the analytic continuation of $I(\lambda)$.

\Theorem IV.12. {\rm (Hilgert-\'Olafsson-\O rsted)} Let 
$(\pi_\lambda, {\cal H}_\lambda)$ be an unitary
highest weight representation of $G$ with $(\pi_\lambda^K, F(\lambda))$ being  
$H\cap K$-spherical. Then $(\pi_\lambda, {\cal H}_\lambda)$
belongs to the relative discrete series of $H\bs G$ if and only if 
the condition 
$$(\forall \alpha\in \Delta_n^+)\quad \la \lambda+\rho, \alpha\ra
<0\leqno({\rm RDS})$$
is satisfied. 

\Proof.  Recall the definition of $\cc^+$, $\a^+$ and the relation 
$C(\b^+)=\cc^+$. Set $A^+\:=\exp_{G^c}(\a^+)$ and let $\|\cdot\|$ denote an arbitrary norm on
$\a$. If we write $(\cc^+)^\star$, then the star $\star$ is to be
taken in $\a^*$.  

\par\nin Step 1: $I(\lambda)<\infty$, if $\lambda+\rho\in
-\Int (\cc^+)^\star$, the interior of $(\cc^+)^\star$. 
\par Here we do not assume that $\lambda\in \a^*$ is dominant integral
with respect to $\Delta_k^+$, but only $\lambda\in \check C_k^\star$. 
By Harish-Chandra's estimates for spherical functions on
non-compact Riemannian symmetric spaces, there exists constants $c>0$ and $d\in
\N$ such that 

$$(\forall \lambda\in \check C_k^\star)
(\forall a\in A^+)\quad \phi_\lambda^0(a)\leq c  a^{\lambda-\rho_k}(1+\|\log a\|)^d\leqno(4.6)$$
(cf.\ [Wal88, 4.5.3]). Note that $J(X)\leq e^{2\rho (C(X))}$ for all
$X\in \b^+$ by the formula for the Jacobian in Lemma IV.6(iii). 
Thus it follows for all $\lambda\in \check C_k^\star$ and
$X=\sum_{j=1}^s x_j X_j \in \b^+$ from (4.6) together with
Lemma IV.7  that 
$$ \eqalign{ \phi_{\lambda+\rho_k}^0(\mu(\exp_G(X))^2) \ J(X) &
\leq c \mu(\exp_G(X))^{2\lambda }(1+ \|\log \mu(\exp_G(X))^2\|)^d e^{2\rho(C(X))}\cr
&\leq c e^{2\lambda(C(X))}(1+2\|C(X)\|)^d e^{2\rho(C(X))}\cr
&\leq c e^{2(\lambda+\rho)(C(X))}(1+2\|C(X)\|)^d.\cr} \leqno(4.7)$$
Now Proposition IV.11 shows that $I(\lambda)<\infty$ if 
$\lambda+\rho \in -\Int (\cc^+)^\star$, proving our first step. 

\par\nin Step 2: $\lambda+\rho\in -\Int (\cc^+)^\star$, if
$I(\lambda)<\infty$. 
\par Recall that $\lambda$ is
supposed to be dominant integral with respect to $\Delta_k^+$. Thus it
follows from (4.4) and the fact that the $H\cap K$-spherical vector $\nu_0$ has a non-zero
$v_\lambda$-component (cf.\ [Hel84, p.\ 537, (7)]) that there is  
a constant $c_\lambda>0$ such that $c_\lambda a^\lambda \leq \phi_{\lambda
+\rho_k}^0(a)$ holds for all $a\in A^+$. Hence  Lemma
IV.7 implies that 

$$(\forall X\in\b^+) \quad {c_\lambda\over 2} e^{2\lambda(C(X))} J(X)\leq
\phi_{\lambda+\rho_k}^0(\mu(\exp_G(X))^2) J(X).$$
In view of Proposition IV.10 and Lemma IV.11, we obtain 
$\lambda+\rho\in -\Int (\cc^+)^\star$ if  
$I(\lambda)<\infty$. This proves our second step. 

\par\nin Step 3: If $\lambda\in \check C_k^\star$, then
$\lambda$ satisfies (RDS) if and only if $\lambda+\rho\in -\Int (\cc^+)^\star$.
\par Note that $\lambda$ satisfies (RDS) means that 
$\lambda+\rho \in -\Int C_{\rm min}^\star$. Now if $\lambda\in \check
C_k^\star$, then $\lambda+\rho \in \Int  \check C_k^\star$. Thus Step 3
follows from Lemma IV.5.  

\par In view of Steps 1-3, it follows that $I(\lambda)$ is finite
if and only if $\lambda$ satisfies the condition (RDS). The proof of
the theorem will therefore be complete with 

\par\nin Step 4: If $\lambda$ satisfies (RDS), then $(\pi_\lambda,
{\cal H}_\lambda)$ is $H$-spherical. 
\par Let $\kappa\: G\to \tilde {K_\C}/ (\tilde {K_\C}\cap \tilde{H_\C})_0$ the
canonical projection defined via the decomposition in Proposition
II.10. Now the function 
$$H\bs G\to \C, \ \ Hg\mapsto \la \pi_\lambda^K(\kappa(g)).v_\lambda, \nu_0\ra$$
generates an $H$-spherical module in the relative discrete series on 
$H\bs G$ (cf.\ [\'O\O91, Th.\ 5.2]). This proves Step 4 and concludes
the proof of the theorem.\qed

The prescription   
$$W\:=-\Int C_{\rm min}^\star\cap \check C_k^\star\subeq -\Int (\cc^+)^\star$$
defines a convex cone in $\a^*$. We write $T_W=i\a^* +W$ for the
associated tube domain in $\a_\C^*$.
Note that $\rho_n\in i\z(\k)^*$ by the construction of $\Delta_n^+$ and
so $-\rho_n\in W$.

\Lemma IV.13.  The function $I(\lambda)$ extends naturally to a
continuous function on the affine subtube $T_W-\rho$, also
denoted by $I$, and which is
holomorphic when restricted to $T_{W^0}-\rho$.  If $m\in \N$ is sufficiently large, then
$W-m\rho_n\subeq W-\rho$ and $I\res_{T_W-m\rho_n}$ is bounded.

\Proof.  First we show that $W-m\rho_n\subeq W-\rho$ for large values
of $m\in \N$. Since $\rho_n\in \Int C_{\rm min}^\star$, we have 
$\rho-m\rho_n\in -\Int C_{\rm min}^\star$ provided $m\in \N$ is
sufficiently large. Further $\rho_n\in i\z(\k)^*$ shows that
$\R.\rho_n\in \check C_k^\star$. Thus we have $\rho-m\rho_n\in W$ if $m$ is chosen
sufficiently large, proving our claim. 
\par Recall the formula for $I(\lambda)$ from Proposition
IV.10. Then (4.7) yields constants $c>0$, $d\in \N$ such that 

$$I(\lambda)\leq
{c\over \dim F(\lambda)}\int_{\cc^+} e^{2(\lambda+\rho)(X)}(1+2\|X\|)^d\ dX\leqno(4.8)$$
holds for some norm $\|\cdot\|$ on $\a$. 
Let $\hat\rho_k$ denote the half sum of the roots in $\hat\Delta_k^+$ and
recall Weyl's Dimension Formula 
$$\dim F(\lambda)={\prod_{\hat\alpha\in \hat\Delta_k^+} \la \lambda+\hat\rho_k,
\hat\alpha\ra \over \prod_{\hat\alpha\in \hat\Delta_k^+} \la \hat\rho_k,
\hat\alpha\ra}.$$
In particular, we see that $\lambda\mapsto {1\over \dim F(\lambda)}$ extends to
a holomorphic map on $T_W$ and $T_W-\rho$ which is bounded when restricted to 
$T_W-m\rho_n$ for all $m\in \N_0$. Further for each fixed $b\in B^+$ the mapping
$$\a_\C^*\to \C, \ \ \lambda\mapsto \phi_{\lambda+\rho_k}^0(\mu(b)^2)$$
is holomorphic. 
Now (4.8) together with Proposition IV.10 imply  that $I(\lambda)$
extends to a continuous function on $T_W-\rho$ which is holomorphic on
$T_{W^0}-\rho$ and bounded when restricted to $T_W-m\rho_n$ provided
$m$ is chosen sufficiently 
large. \qed

\Lemma IV.14. If $m\in \N$ is sufficiently large, then the function 
$$T_{W^0}-m\rho_n\to\C, \ \ \lambda \mapsto c(\lambda+\rho)$$
is holomorphic and bounded. 

\Proof. In view of $\rho_n\in i\z(\k)^*$, this is immediate from Theorem II.14.  \qed

\Theorem IV.15. {\rm (The formal dimension for the relative holomorphic
discrete series on a compactly causal symmetric space)} Let $H\bs G$ be a simply
connected symmetric space associated to a compactly causal symmetric
Lie algebra $(\g,\tau)$ and  $(\pi_\lambda, {\cal H}_\lambda)$
be an unitary highest weight representations of
$G$ for which $F(\lambda)$ is $H\cap K$-spherical. Then the following assertions hold:

\item{(i)} The representation $(\pi_\lambda, {\cal H}_\lambda)$
belongs to the relative discrete series for $H\bs G$ if and only if
the condition 
$$(\forall \alpha\in \Delta_n^+)\quad \la \lambda+\rho, \alpha\ra
<0\leqno({\rm RDS})$$
is satisfied. 
\item{(ii)} If $(\pi_\lambda, {\cal H}_\lambda)$ belongs to the
relative discrete series of $H\bs G$, then the formal dimension
$d(\lambda)$ is given by
$$d(\lambda)=d(\lambda)^G c(\lambda+\rho),$$
where $d(\lambda)^G$ is the formal dimension of $(\pi_\lambda, {\cal H}_\lambda)$
relative to $G$ and $c$ is the
$c$-function of the 
non-compactly $c$-dual space $H^c\bs G^c$ of $H\bs G$ (cf.\ {\rm
Theorem II.14}). Here the right hand side has to be understood as
an analytic continuation of a product of two meromorphic functions.

\Proof. (i) Theorem IV.12.
\par\nin (ii) Let $\hat\rho$ denote the half sum of the elements in
$\hat\Delta^+$ and recall Harish-Chandra's condition for the relative
discrete series on $G$
$$(\forall \hat\alpha\in \hat\Delta_n^+)\quad \la \lambda+\hat\rho, \hat\alpha\ra<0$$ 
(cf.\ [HC56, Lemma 29]) as well as Harish-Chandra's formula for 
the formal dimension
$d(\lambda)^G$ of the relative discrete series on $G$
$$d(\lambda)^G={\prod_{\hat\alpha\in \hat\Delta^+} \la \lambda+\hat\rho,
\hat\alpha\ra \over \prod_{\hat\alpha\in \hat\Delta^+} \la \hat\rho,
\hat\alpha\ra }$$
(cf.\ [HC56, Th.\ 4]). In particular for $m\in \N$ sufficiently 
large, the prescription $\lambda\mapsto
{1\over d(\lambda)^G}$ defines a bounded holomorphic function on 
the affine tube $T_{W^0}-m\rho_n$.

\par Now it follows from Lemma IV.13 and Lemma IV.14 that the 
function 
$$f\: T_{W^0}-m\rho_n  \to \C, \ \ \lambda\mapsto
I(\lambda)c(\lambda+\rho)- {1\over d(\lambda)^G}$$
is holomorphic and bounded for $m$ sufficiently large. 
For such  $m$ Theorem III.6 implies that $f(\lambda)=0$ for all
$\lambda\in W^0-m\rho_n$ which are dominant integral with respect to
$\Delta_k^+$. Thus the identity criterion of Proposition A.2
in Appendix A applies and yields $f=0$. We conclude in particular
that $I(\lambda)^{-1}$ defines a continuation  of
$\lambda\mapsto d(\lambda)^Gc(\lambda+\rho)$ to a continuous function
on  $T_W-\rho$ which is holomorphic when restricted to the interior
$T_{W^0}-\rho$.  Since by definition  $d(\lambda)=I(\lambda)^{-1}$, the assertion
in (ii) follows because $\lambda$ satisfies (RDS) if and only if
$\lambda\in T_W-\rho$.\qed

The following result has already been obtained earlier by Faraut, Hilgert and
\'Olafsson in [FH\'O94, Lemma 9.2], but with a completely different
type of arguments (see  also Theorem II.14).

\Corollary IV.16.  Suppose that $(\g,\tau)=(\h\oplus\h, \sigma)$ is of 
group type (cf.\ {\rm Lemma I.3(i)(2)}). Then the domain of 
convergence ${\cal E}$ for $c$ is given by 

$${\cal E}=i\a^*+(-\Int C_{\rm min}^\star)\cap \Int \check C_k^\star$$
and there exists a constant $\gamma>0$ only depending on the choice 
of the various Haar measures such that 
$$ c(\lambda)=\gamma{1\over \prod_{\alpha\in \Delta^+}
\la \lambda, \alpha\ra}$$
for $\lambda\in {\cal E}$. 

\Proof. In the following we use the notation of Remark III.5. 
Since $(\g,\tau)$ is of group type we have
$d(\lambda)^G=d(\lambda)^{(G_0\times G_0)}=(d(\lambda)^{G_0})^2$, and so it follows from Theorem IV.15(ii)
that $c(\lambda+\rho)={1\over d(\lambda)^{G_0}}$ holds 
for the analytic continuations. In view of Harish-Chandra's formula
for $d(\lambda)^{G_0}$ (cf.\ [HC56, Th.\ 4]), this proves the
corollary.\qed

\Problems. The discrete series on $H\bs G$ are constructed by 
analytic methods, i.e., with generating functions (cf.\ [FJ80], [MaOs84],
[\'O\O91]). But from the algebraic point of view there are still
some interesting open problems. 

\par\nin (a) Using the classification sheme of unitary highest weight 
modules (cf.\ [EHW83]) together with the fine structure theory of compactly 
causal symmetric Lie algebras provided by Theorem IV.4 and [\'Ol91]  one can check case
by case that 
(RDS) implies that $N(\lambda)=L(\lambda)$. In view of Proposition
II.7(ii), this gives a more algebraic proof of the fact that (RDS)
implies that $(\pi_\lambda, {\cal H}_\lambda)$ is $H$-spherical
whenever $(\pi_\lambda^K, F(\lambda))$ is $H\cap K$-spherical. 
The following questions are therefore natural: What is the
algebraic impact of the condition (RDS)? Does there exists 
an analog of the Parthasarathy-condition (cf.\ [EHW83. Prop.\ 3.9])
for the symmetric space setting? 
\par\nin (b) Give a complete classification of $H$-spherical unitary 
highest weight representations. A first step in this direction might be
Proposition II.7(ii) and Remark II.8.\qed

\sectionheadline{V. Applications to holomorphic representation theory}

In this final section we give a second application 
of the Averaging Theorem: We relate the spherical   
character of a spherical unitary highest weight representation of $G$ to
the corresponding spherical functions on the 
$c$-dual space.

\subheadline{Spherical functions and character theory}

\Definition V.1. Let $(\pi_\lambda, {\cal H}_\lambda)$ be an
$H$-spherical unitary highest weight 
representation of $G$. If $0\neq \nu\in ({\cal
H}_\lambda^{-\omega})^H$
and $v_\lambda$ is an highest weight vector, then we define the {\it spherical 
character} $\Theta_\lambda$  of $(\pi_\lambda, {\cal H}_\lambda)$ by 

$$\Theta_\lambda\: S_{\rm max}^0\to\C, \ \ s\mapsto 
{\la v_\lambda, v_\lambda\ra\over |\la \nu, v_\lambda\ra|^2} 
\la \pi_\lambda(s).\nu, \nu\ra.$$ 
Note that $\Theta_\lambda$ is an $H$-biinvariant holomorphic 
function on $S_{\rm max}^0$ (cf.\ [KN\'O97, Lemma V.6]).\qed

\Remark V.2. The particular normalization of $\Theta_\lambda$ has two reasons. First 
that it coincides in the group case (cf.\ Remark III.5) with the 
standard definition, and second because it has the best analytic properties for the 
assignments $\lambda\mapsto \Theta_\lambda(s)$, $s\in S_{\rm max}^0$ (as less poles as 
possible).\qed

\Definition V.3. (Spherical Functions) Recall the definition of the 
domain ${\cal E}_\Omega\subeq \a_\C^*$ (cf.\ Definition II.12). 
If $\lambda\in {\cal E}_\Omega$, then the {\it spherical function 
with parameter} $\lambda$ is defined by 

$$\phi_\lambda\: S_{\rm max}^0\cap HAN\to \C,\ \  s\mapsto
\int_H a_H(sh)^{\lambda-\rho} \ d\mu_H(h)$$
(cf.\ [FH\'O94] or [KN\'O98]). Recall that the defining integrals converge  
absolutely if and only if $\lambda\in {\cal E}_\Omega$ (cf.\ [FH\'O94,
Th.\ 6.3]). \qed

\Theorem V.4. Let $(\pi_\lambda, {\cal H}_\lambda)$ be an 
$H$-spherical unitary highest weight representation of $G$ such that 
$\lambda+\rho\in {\cal E}_\Omega$ holds. Then the spherical character
$\Theta_\lambda$ of $(\pi_\lambda, {\cal H}_\lambda)$ and the
spherical function $\phi_{\lambda+\rho}$ are related by 
$$(\forall s\in S_{\rm max}^0\cap HAN)\quad \Theta_\lambda(s)={1\over
c(\lambda+\rho)} \phi_{\lambda+\rho}(s).$$
In particular, $\phi_{\lambda+\rho}$ extends naturally to a
$H$-biinvariant holomorphic function on $S_{\rm max}^0$. 

\Proof. Since $\lambda+\rho\in {\cal E}_\Omega$, the 
assumption of Theorem II.16 is satisfied and we can rewrite 
$0\neq\nu\in ({\cal H}_\lambda^{-\omega})^H$ as 
$$\nu={\la \nu, v_\lambda\ra  \over \la v_\lambda, v_\lambda\ra c(\lambda +\rho)} 
\int_H \pi_\lambda(h).v_\lambda \ d\mu_H(h).$$
Thus if we replace the first $\nu$ in the definition of
$\Theta_\lambda$ by this expression, we get for all
$s\in S_{\rm max}^0\cap HAN$ that 

$$\eqalign{\Theta_\lambda(s)&={\la v_\lambda, v_\lambda\ra\over |\la \nu, v_\lambda\ra|^2} 
\la \pi_\lambda(s).\nu, \nu\ra\cr 
&={1\over c(\lambda +\rho)}\cdot {1\over \la v_\lambda, \nu\ra}
\int_H\la\pi_\lambda(sh).v_\lambda, \nu\ra\ d\mu_H(h)\cr  
&={1\over c(\lambda +\rho)}\cdot {1\over \la v_\lambda, \nu\ra}
\int_H\la\pi_\lambda(h_H(sh) a_H(sh)n_H(sh)).v_\lambda, \nu\ra\ d\mu_H(h)\cr  
&={1\over c(\lambda +\rho)}\cdot {1\over \la v_\lambda, \nu\ra}
\int_H\la\pi_\lambda(a_H(sh)).v_\lambda, \nu\ra\ d\mu_H(h)\cr  
&={1\over c(\lambda +\rho)}
\int_H a_H(sh)^\lambda d\mu_H(h)\cr  
&={1\over c(\lambda +\rho)}\phi_{\lambda+\rho}(s),\cr}$$
as was to be shown. \qed

\Remark V.5. (a) We remark here that the relation in Theorem V.4 was long time searched
by G. \'Olafsson (cf.\ [\'Ol98, Open Problem 7(1)]). For further 
interesting problems related to this subject we refer to
[Fa98] and [\'Ol98].
\par\nin (b)  The analytic continuation of the relation in Theorem V.4
has been obtained in [HiKr98]. It has 
far reaching consequences for the theory of $G$-invariant 
Hilbert spaces of 
holomorphic functions on $G$-invariant subdomains of the Stein 
manifold $\Xi_{\rm max}^0=G\times_H iW_{\rm max}^0$. In particular, it
implies  the {\it Plancherel Theorem} for these class of Hilbert
spaces (cf.\ [HiKr98]). 
For further information related to this subject we refer to 
[HiKr99b], [KN\'O97], [Kr98,99b] and [Ne99a].\qed

\sectionheadline{Appendix}

\subheadline{A. An identity criterion for bounded analytic functions
on tubes}

\Lemma A.1. Let $\Pi^+\:=\{ z\in \C\: \Im (z)>0\}$ be the upper half
plane and $H^\infty\:=\{ f\in \Hol(\Pi^+)\: \|f\|_\infty<\infty\}$
the Banach space of bounded holomorphic functions on it. Let
$\alpha>0$ and $N=\{ n\alpha i \: n\in \N\}$. Then the
following identity assertion for elements $f$ of $H^\infty(\Pi^+)$
holds: If $f\res_N= 0$, then $f= 0$.

\Proof. Let $D\:=\{ z\in \C\: |z|<1\}$ and $H^\infty(D)=\{ f\in
\Hol(D)\: \|f\|_\infty <\infty\}$. Let $f\in H^\infty(D)$ and
$\{\beta_n\: n\in \N\}$ be subset of zeros of $f$. Then it follows
from [Ru70, Th.\ 15.23] that 

$$f= 0 \quad \hbox{if}\quad \sum_{n=1}^\infty (1-|\beta_n|)=\infty.\leqno(A.1)$$ 
We consider the Cayley transform 

$$c\: \Pi^+\to D,\ \ z\mapsto {z-i \over z+i}, $$
which is an biholomorphic isomorphism, defining an isomorphism of
Banach spaces
$$c_*\: H^\infty(D)\to H^\infty(\Pi^+), \ \ f\mapsto \tilde f=f\circ c.$$
Let $\alpha_n\:=n\alpha i$. Then we have 
$$\beta_n\:=c(\alpha_n)={n\alpha i - i\over n\alpha i + i}=
{n\alpha -1\over n\alpha+1}.$$
Let $N_0\in \N$ be such that $n\alpha -1 >0$ for all $n\geq N_0$. Then

$$\sum_{n=1}^\infty (1-|\beta_n|)\geq \sum_{n=N_0}^\infty (1-{n\alpha
-1\over n\alpha+1}) =\sum_{n=N_0}^\infty {2\over n\alpha+1}
=\infty.\leqno(A.2)$$
Thus if $\tilde f \in H^\infty(\Pi^+)$ vanishes on all $\alpha_n$,
$n\in \N$, then $f(\beta_n)=0$ for all $n\in \N$ and so $f= 0$ by
(A.1) and (A.2). Therefore $\tilde f=c_*(f)=0$, proving the lemma.\qed

\Proposition A.2. Let $\eset\neq W\subeq \R^n$ be an open convex cone,
$T_W\:=\R^n +iW$ the associated tube domain in $\C^n$ and 
$H^\infty(T_W)=\{ f\in \Hol(T_W)\: \|f\|_\infty<\infty\}$
the space of bounded holomorphic functions on $T_W$. Let $\Gamma\subeq
\R^n$ be a lattice. Then the following identity assertion holds:
$$(\forall f\in H^\infty(T_W)) \quad
f\res_{i(\Gamma\cap W)}=0\ \Rarrow\  f= 0.$$

\Proof. We prove the assertion by induction on the dimension $n\in
\N$. 

\par If $n=1$, then $\Gamma=\Z\alpha$ for some $\alpha>0$ and 
$W =\R, \R^+$ or $\R^-$. If  $W=\R$, then the assertion
follows from Liouville's Theorem. In the two remaining cases the
assertion follows from Lemma A.1.

Suppose now the assertion is true for all all dimensions less or equal
to $n-1$, $n\geq 2$. Let $f\in H^\infty(\R^n+iW)$ be an element
vanishing on $i(\Gamma\cap W)$. We have to show that $f= 0$. 
Since $W$ is open, we find a basis $e_1, \ldots, e_n$ of $\R^n$
which is contained in $\Gamma\cap W$. By the Identity Theorem for 
analytic functions, it suffices to prove the assertion for 
$\Gamma=\Z e_1\oplus\ldots\oplus \Z e_n$ and $W=\sum_{j=1}^n \R^+
e_j$. Let $\Gamma_{n-1}=\Z e_1\oplus\ldots\oplus \Z e_{n-1}$ and
$W_{n-1}=\sum_{j=1}^{n-1} \R^+e_j$. Write the variables
$z\in\C^n$ as tuples  $z=(z', z_n)$ with $z'=(z_1, \ldots,
z_{n-1})$. By induction we obtain that $f(z)=f(z', z_n)$ does not
depend on the $z'$-variable. Thus $f(z)=F(z_n)$ for  some $F\in
H^\infty(\Pi^+)$ with $F\res_{\N i}= 0$. Thus by the induction 
hypothesis $F= 0$, and hence $f=0$ establishing the induction step. \qed

\subheadline{B. A lemma on spherical highest weight modules}

Throughout this subsection $(\g,\tau)$ denotes a simple hermitian 
symmetric Lie algebra. Further we use the notation from Section I-II.

\Lemma B.1. Suppose that $(\g,\tau)$ is a simple 
hermitian symmetric Lie algebra and $(G,\tau)$ an associated 
simply connected Lie group. Set $H=G^\tau$ and assume that 
there exist a non-trivial $H$-spherical unitary highest weight representation 
$(\pi_\lambda, {\cal H}_\lambda)$ of $G$. Then the symmetric Lie
algebra $(\g,\tau)$ has to be compactly causal. 

\Proof. Write $\g=\k\oplus \p$ for a $\tau$-invariant Cartan 
decomposition of $\g$  and let $K$ denote the analytic subgroup of $G$ 
corresponding to $\k$. 
\par By assumption we have $({\cal H}_\lambda^{-\omega})^H\neq \{0\}$. 
In particular we can conclude that the module 
$L(\lambda)$ of $K$-finite vectors of $(\pi_\lambda, {\cal H}_\lambda)$
admits nontrivial $H\cap K$-fixed vectors. Recall that $L(\lambda)$ 
is the unique irreducible quotient of the generalized 
Verma module
$$N(\lambda)={\cal U}(\g_\C)\otimes_{{\cal U}(\k_\C\oplus\p^+)} F(\lambda).$$
In particular, there exists an element 
$0\neq v_0\in N(\lambda)^{H\cap K}$. 
\par Recall that $N(\lambda)$ is $\k_\C$-isomorphic to 
${\cal S}(\p^-)\otimes F(\lambda)$, where the $\k_\C$-action 
on ${\cal S}(\p^-)\otimes F(\lambda)$ is defined by 
$$X.(p\otimes v)\:=[X,p]\otimes v+ p\otimes X.v\leqno(B.1)$$
for $X\in\k_\C$, $p\in {\cal S}(\p^-)$ and $v\in F(\lambda)$
(cf.\ [EHW83]).  

\par In order to show that $(\g,\tau)$ is compactly causal, we have 
to prove $\z(\k)\subeq \q$. Assume the contrary, i.e., 
$\z(\k)\subeq \h$. Recall the definition of the element
$Z_0\in\z(\k)$ from Section I and set $X_0\:=-iZ_0\in i\z(\k)$. 
Then the spectrum of $X_0$, considered as an operator 
on the symmetric algebra ${\cal S}(\p^-)$, is $-\N_0$, and we obtain a natural 
grading by homogeneous elements: 
${\cal S}(\p^-)=\bigoplus_{n=0}^\infty {\cal S}(\p^-)^{-n}$. Then 
$N(\lambda)=\bigoplus_{n=0}^\infty {\cal S}(\p^-)^{-n}\otimes F(\lambda)$ 
and we conclude from (B.1) that $X_0$ acts on 
${\cal S}(\p^-)^{-n}\otimes F(\lambda)$ by $-n+\lambda(X_0)$
times the identity. Write $v=\sum_{n=0}^\infty v_0^{-n}$ according to 
the decomposition 
$N(\lambda)=\bigoplus_{n=0}^\infty {\cal S}(\p^-)^{-n}\otimes F(\lambda)$.
Since $X_0\in i(\h\cap\k)$, the element $v_0$ is annihilated by $X_0$ and so 
we must have $v_0=v_0^{-n}$ for some $n\in \N_0$ with 
$\lambda(X_0)=n\geq 0$. But 
a necessary condition for $L(\lambda)$ to be unitarizable is 
$\lambda(X_0)<0$ (cf.\ [Ne99b, Th.\ XI.2.37(ii)]). This gives us 
a contradiction and proves the lemma.\qed

\def\entries{

\[Ba87 Ban, E.\ P. van den, {\it Invariant differential operators on a semisimple symmetric space 
and finite multiplicities in a Plancherel formula}, Ark.\ Mat. {\bf 25} (1987), 175--187

\[BS97 Ban, E.\ P. van den, and  H. Schlichtkrull, {\it The most continuous
part of the Plancherel decomposition for a reductive symmetric space}, 
Ann. of Math. {\bf 145} (1997), 267--364

\[BS99 ---, {\it Fourier inversion on a reductive symmetric space}, 
Acta Math. {\bf 182} (1999), 25--85

\[Be57 Berger, M., {\it Les espaces sym\'etriques non compacts}, Ann.\ Sci.\
\'Ecole Norm.\ Sup. {\bf 74} (1957), 85--177

\[Ch98 Chadli, M., {\it Noyau de Cauchy-Szeg\"o d'un espace
sym\'etrique de type Cayley}, Ann. Inst. Fourier {\bf 48(1)} (1998), 97--132

\[De98 Delorme, P., {\it Formule de Plancherel pour les espaces
sym\'etriques r\'eductifs}, Ann. of Math. {\bf 147} (1998), 417--452

\[EHW83 Enright, T. J., R. Howe, and N. Wallach, {\it A classification of 
unitary highest weight modules} in Proc. ``Representation theory of reductive 
groups" (Park City, UT, 1982), pp. 97-149; Progr. Math. {\bf 40} (1983), 97 
--143 

\[FH\'O94 Faraut, J., J.\ Hilgert, and G. \'Olafsson, {\it Spherical functions
on ordered symmetric spaces}, Ann. Inst. Fourier {\bf 44} (1994), 927--966

\[Fa95 Faraut, J., {\it Fonctions Sph\'eriques sur un Espace
Sym\'etrique Ordonn\'e de Type Cayley}, Contemp. Math. {\bf 191}
(1995), 41--55

\[Fa98 ---, {\it Quelques probl\'emes d'analyse sur les espaces 
sym\'etriques ordonn\'es}, in ``Positivity in Lie Theory: Open Problems'', 
J. Hilgert, J. Lawson, K.--H. Neeb, and E. Vinberg, editors, 
de Gruyter, 1998 

\[FJ80 Flensted-Jensen, M., {\it Discrete series for semisimple
symmetric spaces}, Ann. of Math. {\bf (2)111}(1980), 253--311

\[GiKa62 Gindikin, S. G., and F. I. Karpelevi\v c, {\it Plancherel measure of 
Riemannian symmetric spaces of non-positive 
curvature}, Dokl. Akad. Nauk. SSSR {\bf 145} (1962), 252--255

\[Gr96 Graczyk, P., {\it Espace de Hardy d'un espace sym\'etrique
de type Hermitien}, in ``Journ\'ees Program Gelfand-Gindikin'', Paris,
1996

\[Gr97 ---, {\it Function $c$ on an ordered symmetric space},
Bull. Sci. math. {\bf 121} (1997), 561--572

\[HC56 Harish-Chandra, {\it Representations of semi-simple Lie groups, VI}, 
Amer. J. Math. {\bf 78} (1956), 564--628

\[HC66 ---, {\it Discrete series for semi-simple Lie groups II}, Acta Math.\ {\bf 166}
(1966), 1--111 

\[Hel78 Helgason, S., ``Differential geometry, Lie groups, and symmetric 
spa\-ces,'' Acad. Press, London, 1978

\[Hel84 ---, ``Groups and Geometric Analysis,'' Acad. Press, London,
1984

\[HiKr98 Hilgert, J., and B.\ Kr\"otz, {\it The Plancherel Theorem for 
invariant Hilbert spaces}, submitted 

\[HiKr99a ---,  {\it Representations, characters and spherical 
functions associated to causal symmetric spaces}, J. Funct. Anal., to appear 

\[HiKr99b ---,  {\it Weighted Bergman spaces associated to causal 
symmetric spaces}, manusc. math. {\bf 99 (2)} (1999), 151--180

\[Hi\'Ol96 Hilgert, J.\ and 
G.\ \'Olafsson, ``Causal Symmetric Spaces, Geometry and
Harmonic Analysis,'' Acad. Press, 1996 

\[H\'O\O{}91 Hilgert, J., \'Olafsson, G., and B. \O rsted, {\it Hardy Spaces on 
Affine Symmetric Spaces}, J. reine angew. Math. {\bf 415} (1991), 189--218

\[KoWo65 Kor\'anyi, A., and J. A. Wolf, {\it Realization of hermitean symmetric 
spaces as generalized half planes}, Ann. of. Math. {\bf 81} (1965), 265--288

\[Kr98 Kr\"otz, B., {\it On Hardy and Bergman spaces on complex
Ol'shanski\u\i{} semigroups}, Math. Ann.{\bf 312} (1998), 13-52

\[Kr99a ---, {\it Norm estimates for unitary highest weight
modules}, Ann. Inst. Fourier {\bf 49(4)}, 1--24

\[Kr99b ---, {\it The Plancherel Theorem for Biinvariant Hilbert
Spaces}, Publ.\ RIMS {\bf 35 (1)} (1999), 91--122

\[KrNe96 Kr\"otz, B., and K. - H. Neeb, {\it On hyperbolic cones and
mixed symmetric spaces}, Journal of Lie Theory {\bf 6:1}(1996), 69--146

\[KN\'O97 Kr\"otz, B., K. - H. Neeb, and G. \'Olafsson, {\it  Spherical 
Representations and Mixed Symmetric Spaces}, Represent. Theory {\bf
1} (1997), 424-461

\[KN\'O98 ---, {\it Spherical functions on mixed symmetric spaces }, submitted

\[Kr\'Ol99 Kr\"otz, B., and G.\ \'Olafsson, {\it The c-function for non-compactly 
causal symmetric spaces}, submitted 

\[La94 Lawson, J.D., {\it Polar and Ol'shanski\u\i{} decompositions}, J. f\"ur 
Reine Ang. Math. {\bf 448} (1994), 191--219

\[Lo69 Loos, O., ``Symmetric Spaces I : General Theory", Benjamin, New York, 
Amsterdam, 1969

\[MaOs84 Matsuki, T., and T. Oshima, {\it A description of discrete
series for semisimple symmetric spaces}, Adv. Stud. Pure Math. {\bf
4}, 1984, 229-390 

\[Mo64 Moore, C.C., {\it Compactifications of symmetric spaces, II,
The Cartan domains}, Amer. J. Math. {\bf 86} (1964), 358--378

\[Ne99a Neeb, K.--H., {\it On the complex geometry of invariant domains in
complexified symmetric spaces}, Ann. Inst. Fourier {\bf 49(1)} (1999), 
177--225

\[Ne99b  ---,  ``Holomorphy and Convexity in Lie Theory,'' 
Expositions in Mathematics, de Gruyter, in press

\[\'Ol87 \'Olafsson, G., {\it Fourier and Poisson transformation associated to a
semsisimple symmetric space}, Invent. math. {\bf 90} (1987), 605--629

\[\'Ol91 ---, {\it Symmetric Spaces of Hermitean Type}, Differential 
Geometry and its Applications {\bf 1} (1991),195--233

\[\'Ol97 ---, {\it Spherical Functions and Spherical 
Laplace Transform on Ordered Symmetric Space}, submitted 

\[\'Ol98 ---, {\it Open Problems in Harmonic Analysis on Causal 
Symmetric Spaces}, in ``Positivity in Lie Theory: Open Problems'', 
J. Hilgert, J. Lawson, K.--H. Neeb, and E. Vinberg, editors, 
de Gruyter, 1998 

\[\'O\O91 \'Olafsson, G., and B. \O{}rsted, {\it The holomorphic discrete series of  
affine symmetric spaces and representations with reproducing kernels},
Trans. Amer. Math, Soc. {\bf 326} (1991), 385-405

\[Ru70 Rudin, W., ``Real and Complex Analysis'', McGraw Hill, London, New
York, 1970

\[Sch84 Schlichtkrull, H., ``Hyperfunctions and Harmonic Analysis on Symmetric\break 
Spaces", Progress in Math. 84, Birkh\"auser, 1984 

\[Wal88 Wallach, N., ``Real Reductive Groups I,'' Academic Press, 1988
}

{\sectionheadline{\bf References}
\frenchspacing
\entries\par}
\dlastpage
\vfill\eject
\bye